\documentclass[a4paper,11pt,twoside]{article}
\usepackage{amssymb,amsmath,latexsym,geometry,fancyhdr,lineno,hyperref,titletoc,caption}
\usepackage{tikz}
\usepackage{float}
\tikzset{elegant/.style={smooth,thick,samples=50,black}}
\tikzset{eaxis/.style={->,>=stealth}}

\contentsmargin{0pt}

\dottedcontents{section}[2em]{\vspace{-1mm}\small}{2em}{0pt}
\dottedcontents{subsection}[5em]{\vspace{-1mm}\small}{3em}{5pt}

\hypersetup{colorlinks,linkcolor= blue,citecolor=blue}
\newtheorem{theorem}{Theorem}[section]

\newtheorem{lemma}[theorem]{Lemma}
\newtheorem{proposition}[theorem]{Proposition}
\newtheorem{remark}{Remark}[section]

\newtheorem{definition}[theorem]{Definition}
\def\proof{\mbox {\textbf{Proof.}~~}}
\numberwithin{equation}{section}
\allowdisplaybreaks[0]

\begin{document}
\title{{\bf\Large Normalized solutions with positive energies for a coercive problem and application to the cubic-quintic nonlinear  Schr\"{o}dinger equation}}
\author{\\
{ \textbf{\normalsize Louis Jeanjean}}\\%\footnote{}\\%\footnote{Corresponding author}
{\it\small Laboratoire de Math\'{e}matiques (CNRS UMR 6623),}\\
{\it\small Universit\'{e} de Bourgogne Franche-Comt\'{e},}\\
{\it\small Besan\c{c}on 25030, France}\\
{\it\small e-mail: louis.jeanjean@univ-fcomte.fr}\\
\\
{ \textbf{\normalsize Sheng-Sen Lu}}\\%\footnote{}\\%\footnote{Corresponding author}
{\it\small LMAM and School of Mathematical Sciences,}\\
{\it\small Peking University,}\\
{\it\small Beijing 100871, PR China}\\
{\it\small e-mail: sslu@pku.edu.cn}}
\date{}
\maketitle
{\bf\normalsize Abstract.}
{\small
  In any dimension $N \geq 1$, for given mass $m > 0$ and when the $C^1$ energy functional
    \begin{equation*}
      I(u) := \frac{1}{2} \int_{\mathbb{R}^N} |\nabla u|^2 dx - \int_{\mathbb{R}^N} F(u) dx
    \end{equation*}
  is coercive on the mass constraint
    \begin{equation*}
      S_m := \left\{ u \in H^1(\mathbb{R}^N) ~|~ \|u\|^2_{L^2(\mathbb{R}^N)} = m \right\},
    \end{equation*}
  we are interested in searching for constrained critical points at positive energy levels. Under general conditions on $F \in C^1(\mathbb{R}, \mathbb{R})$ and for suitable ranges of the mass, we manage to construct such critical points which appear as a local minimizer or correspond to a mountain pass or a symmetric mountain pass level. In particular, our results shed some light on the cubic-quintic nonlinear Schr\"{o}dinger equation in $\mathbb{R}^3$.
}



%%\let\thefootnote\relax\footnotetext{{\bf \footnotesize Date:} { \footnotesize \date{\today}}}
%%\let\thefootnote\relax\footnotetext{{\bf \footnotesize 2010 MSC:} { \footnotesize }}
%%\let\thefootnote\relax\footnotetext{{\bf \footnotesize Key words:} { \footnotesize }}

%%%%----------------------------------------------------------------------------------------------------------------------------
%%%%----------------------------------------------------------------------------------------------------------------------------
\pagestyle{fancy}
\fancyhead{} % clear all header fields
\fancyfoot{} % clear all footer fields
\renewcommand{\headrulewidth}{0pt}
\renewcommand{\footrulewidth}{0pt}
\fancyhead[CE]{ \textsc{Louis Jeanjean \& Sheng-Sen Lu}}
\fancyhead[CO]{ \textsc{Normalized Positive Energy Solutions}}
\fancyfoot[C]{\thepage}
%%%%----------------------------------------------------------------------------------------------------------------------------
%%%%----------------------------------------------------------------------------------------------------------------------------

%%\tableofcontents
%%\pagewiselinenumbers
%%%%----------------------------------------------------------------------------------------------------------------------------
%%%%----------------------------------------------------------------------------------------------------------------------------

\section{Introduction}

This paper is concerned with the nonlinear scalar field equation with an $L^2$ constraint
  \begin{equation*}\tag{$P_m$}\label{problem}
    \left\{
      \begin{aligned}
        - \Delta u & = f(u) - \mu u \quad \text{in} ~ \mathbb{R}^N,\\
        \|u\|^2_{L^2(\mathbb{R}^N)} & = m,\\
        u & \in H^1(\mathbb{R}^N),
      \end{aligned}
    \right.
  \end{equation*}
where $N \geq 1$, $f \in C(\mathbb{R}, \mathbb{R})$, $m > 0$ is a given mass, $H^1(\mathbb{R}^N)$ is the Sobolev space, and $\mu \in \mathbb{R}$ is not fixed  {\it a priori} but appears as an unknown Lagrange multiplier due to the mass constraint $\|u\|^2_{L^2(\mathbb{R}^N)} = m$. A function $u \in H^1(\mathbb{R}^N)$ is said to be a solution of \eqref{problem} if for some $\mu \in \mathbb{R}$ the couple $(u, \mu)$ satisfies \eqref{problem}.

The main feature of \eqref{problem} is that the desired solutions have an a priori prescribed $L^2$-norm. Solutions of this type are often referred to in the literature as normalized solutions. As may be well known, the study of \eqref{problem} naturally arises in the search of standing waves for the time dependent nonlinear Schr\"{o}dinger equation
  \begin{equation}\label{eq:equation-evolution}
    i \partial_t \psi + \Delta \psi + f(\psi) =0, \qquad \psi : \mathbb{R}  \times \mathbb{R}^N \to \mathbb{C}.
  \end{equation}
By standing waves, we mean solutions to \eqref{eq:equation-evolution} of the special form $\psi(t,x) := e^{i\mu t} u(x)$ with $\mu \in \mathbb{R}$ and $u \in H^1(\mathbb{R}^N)$. Clearly, the function $\psi(t,x) := e^{i\mu t} u(x)$ solves \eqref{eq:equation-evolution} if and only if the profile $u(x)$ satisfies
  \begin{equation}\label{eq:equation-libre}
    - \Delta u = f(u) - \mu u  \quad \mbox{in} ~ H^1(\mathbb{R}^N)
  \end{equation}
for the associated real number $\mu  \in \mathbb{R}$. The investigation of such type of equations, which had already a strong motivation forty years ago,  see for instance \cite{Be83-1,Be83-2,Lions84-1,Lions84-2,St77},  now lies at the root of several models directly linked with current applications (such as nonlinear optics, the theory of water waves, and Bose-Einstein condensates). For these equations, finding solutions with a prescribed $L^2$-norm is particularly relevant since this quantity is preserved along the time evolution.

The problem \eqref{problem} is variational in nature and thus our basic strategy to study it is to make use of critical points theory. Indeed, under mild conditions on the nonlinearity $f$, one can introduce the $C^1$ energy functional
  \begin{equation*}
    I(u) := \frac{1}{2} \int_{\mathbb{R}^N} |\nabla u|^2 dx - \int_{\mathbb{R}^N} F(u) dx
  \end{equation*}
on $H^1(\mathbb{R}^N)$, where $F(t) := \int^t_0 f(\tau) d\tau$ for any $t \in \mathbb{R}$. For given mass $m > 0$, let
  \begin{equation*}
    S_m := \left\{u \in H^1(\mathbb{R}^N) ~ | ~ \|u\|^2_{L^2(\mathbb{R}^N)} = m \right\}.
  \end{equation*}
It is easily seen that solutions to \eqref{problem} can be obtained as critical points of the functional $I$ constrained to the sphere $S_m$. For future reference, following \cite{BJ16} we say that $v \in S_m$ is an energy ground state solution to \eqref{problem} if it is a solution of \eqref{problem} that minimizes the value of $I$ among all the solutions of \eqref{problem}. Namely, if
  \begin{equation*}
    \big(I_{|S_m}\big)'(v) = 0 \qquad \text{and} \qquad I(v) =  \inf \big\{I(u) ~|~ \big(I_{|S_m}\big)'(u) = 0 \big\}.
  \end{equation*}
Note that this definition keeps its meaning even if the global infimum
  \begin{equation}\tag{$\text{Inf}_m$}\label{eq:Inf_m}
    E_m := \inf_{u \in S_m} I(u)
  \end{equation}
is finite but not achieved or the constrained functional $I_{|S_m}$ is unbounded from below.

In the present paper we shall revisit the situation when the constrained functional $I_{|S_m}$ is bounded from below and coercive for any $m > 0$. This case is commonly called mass subcritical and has been under extensive studies for decades. Among many possible choices,  we refer the reader to \cite{CS21, C03, CL82, HT19, JL19, JL21, LRN20, Sh14,  St19}  and to the references therein. In contrast however, to these previous works, we are herein interested in searching for constrained critical points at positive energy levels. In fact, for a coercive functional on an $L^2$ constraint, associated with an equation such as \eqref{eq:equation-libre}, this seemingly simple issue has not so far been much considered. Compared with the common case where one looks for constrained critical points having negative energies, one faces several new difficulties in the positive energy scenario. First, one needs to identify where, on the $L^2$ constraint, to look for such critical points. Most of the time this step relies on a good understanding of the geometry of the constrained functional. This being done, one is led to study the convergence of particular sequences, typically minimizing or Palais-Smale sequences.  Although these sequences must be bounded due to the fact that the constrained functional is coercive, it often turns out challenging to obtain their strong convergence. In particular, showing that the vanishing may not occur is delicate at positive energy levels while at negative ones this step is most of the time straightforward.

To clarify our mass subcritical setting and for convenience of later discussion, we give a brief account of our recent work \cite{JL21}. Among other issues like the least action characterization and the radial symmetry of global minimizers, the solvability of the global minimization problem \eqref{eq:Inf_m} was investigated there under the following assumptions.
  \begin{itemize}
    \item[$(f1)$] \hypertarget{f1}{} $\lim_{t \rightarrow 0} f(t)/t = 0$.
	\item[$(f2)$] \hypertarget{f2}{} When $N \geq 3$,
                    \begin{equation*}
                      \limsup_{t \to \infty} \frac{|f(t)|}{|t|^{\frac{N + 2}{N - 2}}} < \infty;
                    \end{equation*}
                  when $N = 2$,
                    \begin{equation*}
                      \lim_{t \to \infty} \frac{f(t)}{e^{\alpha t^2}} = 0 \quad \forall \alpha > 0;
                    \end{equation*}
                  and also for any $N \geq 1$,
                    \begin{equation*}
                      \limsup_{t \to \infty} \frac{f(t)t}{|t|^{2 + \frac{4}{N}}} \leq 0.
                    \end{equation*}
    \item[$(f3)$] \hypertarget{f3}{} There exists $\zeta \neq 0$ such that $F(\zeta) > 0$.
  \end{itemize}

\begin{theorem}[{\cite[Theorem 1.1]{JL21}}]\label{th-1}
  Assume that $N \geq 1$ and $f \in C(\mathbb{R}, \mathbb{R})$ satisfies $(\hyperlink{f1}{f1}) - (\hyperlink{f3}{f3})$. Then
    \begin{equation*}
      E_m := \inf_{u \in S_m} I(u) > - \infty
    \end{equation*}
  and the map $m \mapsto E_m$ is nonincreasing and continuous. Moreover,
	\begin{itemize}
      \item[$(i)$] there exists a number $m^* \in [0, \infty)$ such that
                     \begin{equation*}
                       E_m = 0 \quad \text{if} ~ 0 < m \leq m^*, \qquad E_m < 0 \quad \text{when} ~ m > m^*;
                     \end{equation*}
      \item[$(ii)$] when $m > m^*$, the global infimum $E_m$ is achieved and thus \eqref{problem} has an energy ground state solution $v \in S_m$ with $I(v) = E_m < 0$;
	  \item[$(iii)$] when $0 < m < m^*$, $E_m = 0$ is not achieved;
	  \item[$(iv)$] $m^* = 0$ if in addition
                      \begin{equation*}\tag{A.1}\label{eq:key_1}
                        \lim_{t \to 0} \frac{F(t)}{|t|^{2 + \frac{4}{N}}} = + \infty,
                      \end{equation*}
                    and $m^* > 0$ if in addition
                      \begin{equation*}\tag{A.2}\label{eq:key_2}
                        \limsup_{t \to 0}\frac{F(t)}{|t|^{2 + \frac{4}{N}}} < + \infty.
                      \end{equation*}
    \end{itemize}
\end{theorem}

In this paper, we pursue the study of \eqref{problem} within the above general mass subcritical framework and aim to reveal the rich structure of the set of positive energy solutions under the following two additional but simple conditions:
  \begin{itemize}
    \item[$(f4)$] \hypertarget{f4}{} $\limsup_{t \to 0} F(t)/|t|^{2 + 4/N} \leq 0$,
    \item[$(f5)$] \hypertarget{f5}{} there exists $\vartheta \in (2, 2^*)$ such that $\vartheta F(t) \geq f(t)t$ for any $t\in\mathbb{R}$.
  \end{itemize}
Here $2^* := 2N/(N - 2)$ when $N \geq 3$ and $2^* := + \infty$ if $N = 1, 2$. As one may see, the condition $(\hyperlink{f4}{f4})$ shows that $F$ is mass supercritical at zero and thus the constrained functional $I_{|S_m}$ satisfies the estimate in Lemma \ref{lemma:I} $(iii)$. This estimate, which has a simple  geometric interpretation, is essential for us to identify suspected positive critical levels. As to the condition $(\hyperlink{f5}{f5})$, it comes into use only when we search for positive energy solutions of saddle type and will be exploited in our compactness arguments to ensure that the relevant Lagrange multipliers are positive; see Lemma \ref{lemma:PSP_c} and Remark \ref{rmk:PSP_c}. We postpone to Remark \ref{rmk:examples} some examples of nonlinearities satisfying $(\hyperlink{f1}{f1})-(\hyperlink{f5}{f5})$. Since $(\hyperlink{f4}{f4})$ implies \eqref{eq:key_2}, let us keep in mind from now on that $m^* > 0$.

Our first main result shows that under the assumptions $(\hyperlink{f1}{f1})-(\hyperlink{f4}{f4})$ the problem \eqref{problem} admits a positive energy solution in a suitable left neighborhood of $m^* > 0$, even though the global infimum $E_m = 0$ is not achieved for any $m \in (0, m^*)$. It is obtained as a local minimizer, has a constant sign, and in particular proves under an extra proper condition to be an energy ground state solution.
\begin{theorem}\label{theorem:local}
  Assume that $N \geq 1$ and $f \in C(\mathbb{R}, \mathbb{R})$ satisfies $(\hyperlink{f1}{f1})-(\hyperlink{f4}{f4})$. Let $\rho = \rho(m^*) > 0$ be the number given in Lemma \ref{lemma:I} $(iii)$. For any $m \in (0, m^*]$, we set
    \begin{equation*}
      S^\rho_m := \big\{ u \in S_m~|~\|\nabla u\|^2_{L^2(\mathbb{R}^N)} > \rho(m^*) \big\}
    \end{equation*}
  and define the local infimum
    \begin{equation*}
      \overline{E}_m := \inf_{u \in S^\rho_m} I(u) \geq E_m = 0.
    \end{equation*}
  Then the map $m \mapsto \overline{E}_m$ is continuous and nonincreasing. Moreover there exists $m^{**} \in (0, m^*)$ such that when $m \in (m^{**},m^*]$ the local infimum $\overline{E}_m$ is achieved by some $v \in S^\rho_m$ with
    \begin{equation*}
      I(v)=
        \left\{
          \begin{aligned}
          & ~ \overline{E}_m > 0, \qquad \text{for} ~ m \in (m^{**}, m^*),\\
          & ~ \overline{E}_m = 0, \qquad \text{for} ~ m = m^*.
        \end{aligned}
      \right.
    \end{equation*}
  In particular,  $v \in S^\rho_m$ has a constant sign on $\mathbb{R}^N$ and is  a solution of \eqref{problem} with the associated Lagrange multiplier being positive, and $\overline{E}_m$ is strictly decreasing in $m \in (m^{**}, m^*]$. If in addition
    \begin{equation*}\tag{A.3}\label{eq:key_3}
      \limsup_{t \to 0} \frac{f(t)t - 2F(t)}{|t|^{2+4/N}} \leq 0,
    \end{equation*}
  then the local minimizer $v \in S^\rho_m$ is an energy ground state solution of \eqref{problem}.
\end{theorem}

\begin{remark}\label{rmk:comments}
  \begin{itemize}
    \item[$(i)$] As it will be clear from the proof of Theorem \ref{theorem:local},  when $m \in (m^{**}, m^*]$  it also holds that any minimizing sequence for $\overline{E}_m$ is, up to a subsequence and up to translations in $\mathbb{R}^N$, strongly convergent.
    \item[$(ii)$] It is cumbersome, under our general assumptions, to try to derive a sharp estimate on the minimal value $m^{**} >0$ for the existence of minimizers to $\overline{E}_m$. However, a minimal mass threshold appears naturally even for solutions to \eqref{problem}, see Remark  \ref{rmk:no small solution}.
    \item[$(iii)$] The existence of critical points at positive energy levels,  in particular of local minimizers type, was seldom observed in mass subcritical problems. We refer to \cite{JLW15} for a result on a quasilinear problem and to \cite{SS21} which deals with a Schr\"{o}dinger-Poisson equation.
    \item[$(iv)$]  The condition $(\hyperlink{f4}{f4})$ is necessary in the sense that just under $(\hyperlink{f1}{f1})-(\hyperlink{f3}{f3})$ the problem \eqref{problem} may not have positive energy solutions.  Remark \ref{rmk:2D case} provides such a counterexample.
  \end{itemize}
\end{remark}

When $N \geq 2$ and supposing additionally $(\hyperlink{f5}{f5})$, our next result establishes the existence of a positive energy solution of saddle type for any $m > m^{**}$. To avoid misunderstanding when looking at Theorem \ref{theorem:MP} below, one should note that $(\hyperlink{f4}{f4})$ and $(\hyperlink{f5}{f5})$ imply \eqref{eq:key_3}.

\begin{theorem}\label{theorem:MP}
  Assume that $N \geq 2$ and $f \in C(\mathbb{R}, \mathbb{R})$ satisfies $(\hyperlink{f1}{f1})-(\hyperlink{f5}{f5})$. Then for any $m > m^{**}$ the problem \eqref{problem} admits a radial solution $w \in S_m$ which corresponds to a mountain pass level and satisfies
    \begin{equation*}
      \left\{
             \begin{aligned}
               & I(w) > 0 > I(v), & \qquad  & \text{for} ~ m > m^*, \\
               & I(w) > 0 = I(v), & \qquad & \text{for} ~ m = m^*,\\
               & I(w) > I(v) > 0, & \qquad & \text{for} ~ m \in (m^{**}, m^*).
             \end{aligned}
      \right.
    \end{equation*}
  Here $v \in S_m$ is the ground state solution of \eqref{problem} given by Theorems \ref{th-1} and \ref{theorem:local} when $m > m^*$ and $m \in (m^{**}, m^*]$ respectively. If $f$ is also odd, then $w \in S_m$ can be assumed to be nonnegative.
\end{theorem}

Under the setting of Theorem \ref{theorem:MP} and when $f$ is odd, we are able to prove further the following Theorem \ref{theorem:multiplicity}. It seems to be the first existence result of arbitrarily finitely many positive energy solutions in the investigation of $L^2$ constrained mass subcritical problems. In particular, it complements the recent related works \cite{HT19,JL19} where multiple radial solutions of \eqref{problem} were obtained at negative energy levels.
\begin{theorem}\label{theorem:multiplicity}
  Assume that $N \geq 2$ and $f \in C(\mathbb{R}, \mathbb{R})$ is odd satisfying $(\hyperlink{f1}{f1})-(\hyperlink{f5}{f5})$. Then for each $k \in \mathbb{N}^+$ there exists $m_k>0$ such that when $m > m_k$ the problem \eqref{problem} has at least $k$ distinct radial solutions with positive energies.
\end{theorem}

\begin{remark}\label{rmk:examples}
  Let us now give some examples of nonlinearities satisfying our assumptions.
    \begin{itemize}
      \item[$(i)$] $f(t) = a |t|^{p-2}t -b |t|^{q-2}t$ with $a>0, b>0$ and
                     \begin{equation*}
                       \left\{
                         \begin{aligned}
                           & 2 + \frac{4}{N} < p < q < \infty, & \qquad  & \text{if} ~ N = 1,2,\\
                           & 2 + \frac{4}{N} < p < q \leq \frac{2N}{N-2}, & \qquad & \text{if} ~ N \geq 3,
                         \end{aligned}
                       \right.
                     \end{equation*}
                   fulfills $(\hyperlink{f1}{f1}) - (\hyperlink{f5}{f5})$ and \eqref{eq:key_3}.
	  \item[$(ii)$] $f(t) = a |t|^{p - 2}t - b |t|^{q - 2}t - c |t|^{r - 2}t$ with $a > 0, b > 0, c > 0$ and
                      \begin{equation*}
                        \left\{
                          \begin{aligned}
                            & 2 < r < p < q < \infty, & \qquad  & \text{if} ~ N = 1,2,\\
                            & 2 < r < p < q \leq \frac{2N}{N-2}, & \qquad & \text{if} ~ N \geq 3,
                          \end{aligned}
                        \right.
                      \end{equation*}
                    satisfies $(\hyperlink{f1}{f1}), (\hyperlink{f2}{f2}), (\hyperlink{f4}{f4})$ and \eqref{eq:key_3}. Also, it is clear that $(\hyperlink{f3}{f3})$ holds if and only if
                      \begin{equation*}
                        \left[\frac{a}{p( q - r)}\right]^{q - r} >  \left[\frac{b}{q(p - r)}\right]^{p - r} \left[\frac{c}{r(q - p)}\right]^{q - p}.
                      \end{equation*}	
    \end{itemize}
  The above examples are only some special cases, our main theorems apply to more general nonlinearities, in particular to those which are not a sum of competing powers.
\end{remark}

Finally, as an interesting application, we consider \eqref{problem} with the specific choice $N = 3$ and $f(t) = |t|^2 t- |t|^4 t$.  In this case, \eqref{eq:equation-evolution} is called the three-dimensional cubic-quintic nonlinear Schr\"{o}dinger equation. This equation which appears in several physical models, see \cite{M19} for precise references, has been the object of intensive studies these last years, see in particular \cite{CKS21,CS21, KMV21, KOPV17, LRN20, TVZ07}. It is somehow surprising that although  our theorems are obtained under general assumptions on the nonlinearity $f$ they prove useful to derive new results on this specific problem which, as we shall recall in Section \ref{sect:application}, enjoys several convenient properties. See Theorem \ref{theorem:cubic-quintic} for our results in that direction.

The paper is organized as follows. In Section \ref{sect:preliminaries} we establish some preliminary results. In Section \ref{sect:local} we prove Theorem \ref{theorem:local} combining, in a way we have not encountered before, local minimization techniques and robust compactness arguments. Sections \ref{sect:MP} and \ref{sect:SMP} are devoted to the proofs of Theorems \ref{theorem:MP} and \ref{theorem:multiplicity} respectively. This is done through the development of appropriate mountain pass or symmetric mountain pass arguments. Finally, we present in Section \ref{sect:application} some consequences of the above results to the three-dimensional cubic-quintic nonlinear Schr\"{o}dinger equation.

%%%%%%%%%%%%%%%%%%%%%%%%%%%%%%%%%%%%%%%%%%%%%%%%%%%%%%%%%%%%%%%%%%%%%%%%%%%%%%%%%%%%%%%%%%%%%%%%%%%%%%%%%%%%%%%%%%%%%%%%%%%%%%%%%%%%%%
%%%%%%%%%%%%%%%%%%%%%%%%%%%%%%%%%%%%%%%%%%%%%%%%%%%%%%%%%%%%%%%%%%%%%%%%%%%%%%%%%%%%%%%%%%%%%%%%%%%%%%%%%%%%%%%%%%%%%%%%%%%%%%%%%%%%%%

\section{Preliminary results}\label{sect:preliminaries}

In preparation of the proofs of our main theorems, this section contains several technical results.
\begin{lemma}
  Assume that $N \geq 1$ and $f \in C(\mathbb{R}, \mathbb{R})$ satisfies $(\hyperlink{f1}{f1})-(\hyperlink{f2}{f2})$. Then any solution  $u \in S_m$ to \eqref{problem} satisfies $P(u) = 0$ where
    \begin{equation*}
      P(u) := \int_{\mathbb{R}^N} |\nabla u|^2 dx - \frac{N}{2} \int_{\mathbb{R}^N} \big[ f(u)u - 2F(u) \big] dx.
    \end{equation*}
\end{lemma}
\proof The Pohozaev identity corresponding to
  \begin{equation}\label{eq:libre}
    - \Delta u = f(u) - \mu u
  \end{equation}
reads as follows, see \cite[Proposition 1]{Be83-1},
  \begin{equation}\label{eq:pohozaev-libre}
    Q(u) := \frac{N-2}{2N} \int_{\mathbb{R}^N}|\nabla u|^2dx + \frac{\mu}{2} \int_{\mathbb{R}^N}|u|^2dx - \int_{\mathbb{R}^N} F(u) dx = 0.
  \end{equation}
Also, multiplying \eqref{eq:libre} by $u$ and  integrating over $\mathbb{R}^N$  we get
  \begin{equation}\label{eq:Nehari}
    \int_{\mathbb{R}^N} |\nabla u|^2 dx  =  \int_{\mathbb{R}^N}  f(u)u  dx  - \mu \int_{\mathbb{R}^N}  |u|^2  dx.
  \end{equation}
Now, combining \eqref{eq:pohozaev-libre} and \eqref{eq:Nehari}  we obtain  that $P(u) =0$. \hfill $\square$

Our next lemma gives some estimates, whose proofs are standard and left to the reader.
\begin{lemma}\label{lemma:I}
  Assume that $N \geq 1$ and $f \in C(\mathbb{R}, \mathbb{R})$ satisfies $(\hyperlink{f1}{f1})-(\hyperlink{f2}{f2})$. Then the following statements hold.
    \begin{itemize}
	  \item[$(i)$] For any bounded sequence $\{u_n\}$ in $H^1(\mathbb{R}^N)$,
                     \begin{equation*}
                       \lim_{n \to \infty} \int_{\mathbb{R}^N} F(u_n) dx = 0
                     \end{equation*}
                   if $\lim_{n \to \infty} \|\nabla u_n\|_{L^2(\mathbb{R}^N)} = 0$, and
                     \begin{equation*}
                       \limsup_{n \to \infty} \int_{\mathbb{R}^N} F(u_n) dx \leq 0
                     \end{equation*}
                   if $\lim_{n \to \infty} \|u_n\|_{L^{2 + 4/N}(\mathbb{R}^N)} = 0$.
      \item[$(ii)$] There exists $C = C(f, N, m) > 0$  such that
                      \begin{equation*}
                        I(u) \geq \frac{1}{4}\|\nabla u\|^2_{L^2(\mathbb{R}^N)} - C(f, N, m)
                      \end{equation*}
                    for any $u \in H^1(\mathbb{R}^N)$ satisfying $\|u\|^2_{L^2(\mathbb{R}^N)} \leq m$. In particular, $I$ is coercive on $S_m$.
      \item[$(iii)$] For any $m > 0$ there exists $\rho = \rho(m) > 0$ small enough such that,  for all $u \in H^1(\mathbb{R}^N)$ satisfying both $\|u\|^2_{L^2(\mathbb{R}^N)} \leq m$ and $\|\nabla u\|^2_{L^2(\mathbb{R}^N)} \leq 4\rho(m)$, we have respectively
            \begin{equation*}
              I(u) \geq \frac{1}{4} \|\nabla u\|^2_{L^2(\mathbb{R}^N)} \qquad \text{if} ~ (\hyperlink{f4}{f4}) ~ \text{is satisfied},
            \end{equation*}
          and
            \begin{equation*}
              P(u) \geq \frac{1}{2} \|\nabla u\|^2_{L^2(\mathbb{R}^N)} \qquad \text{when} ~ \eqref{eq:key_3} ~ \text{holds}.
            \end{equation*}
    \end{itemize}
  \end{lemma}

The following result answers positively the question of existence of global minimizers for the sharp threshold mass $m^* > 0$ which was left open in Theorem \ref{th-1}. It will not only play an essential role in our proofs of Theorems \ref{theorem:local} and \ref{theorem:MP} but also is interesting by itself.
\begin{lemma}\label{lemma:m*}
  Assume that $N \geq 1$ and $f \in C(\mathbb{R}, \mathbb{R})$ satisfies $(\hyperlink{f1}{f1})-(\hyperlink{f4}{f4})$. Then the global infimum $E_{m^*} = 0$ is  achieved  by some $v \in S_{m^*}$. In particular, the minimizer $v \in S_{m^*}$ is an energy ground state solution  of \eqref{problem} when $m = m^*$ with the associated Lagrange multiplier being positive.
\end{lemma}
\proof Let $m_k := m^* + 1/k$ for any $k \in \mathbb{N}^+$. Since $E_{m_k} < 0$ by Theorem \ref{th-1} $(i)$,  one may choose $v_k \in S_{m_k}$ such that
  \begin{equation}\label{eq:v_k}
    E_{m_k} \leq I(v_k) \leq \frac{1}{2} E_{m_k} < 0 \qquad \text{for each}~ k \in \mathbb{N}^+.
  \end{equation}
Hence the sequence $\{v_k\}$ is bounded in $H^1(\mathbb{R}^N)$ via Lemma \ref{lemma:I} $(ii)$, and up to a subsequence $\lim_{k \to \infty}\int_{\mathbb{R}^N}|\nabla v_k|^2 dx$ and  $\lim_{k \to \infty} \int_{\mathbb{R}^N} F(v_k) dx$ exist. In view of \eqref{eq:v_k} and Lemma \ref{lemma:I} $(iii)$, it is also clear that
  \begin{equation}\label{eq:nalba-vk}
    \lim_{k \to \infty} \|\nabla v_k\|_{L^2(\mathbb{R}^N)} > 0.
  \end{equation}
We prove next that $\{v_k\}$ is non-vanishing, that is
  \begin{equation*}
    \sigma := \limsup_{k \to \infty} \left(\sup_{y \in \mathbb{R}^N} \int_{B(y, 1)}|v_k|^2 dx \right) > 0.
  \end{equation*}
Otherwise, $\sigma = 0$ and then \cite[Lemma I.1]{Lions84-2} gives that $\|v_k\|_{L^{2 + 4/N}(\mathbb{R}^N)} \to 0$. By Lemma \ref{lemma:I} $(i)$ and \eqref{eq:nalba-vk}, we obtain
  \begin{equation*}
    \lim_{k \to \infty} I(v_k) \geq \frac{1}{2}\lim_{k \to \infty}\|\nabla v_k\|^2_{L^2(\mathbb{R}^N)} > 0.
  \end{equation*}
This contradicts \eqref{eq:v_k} and hence $\sigma > 0$. Since $\{v_k\}$ is non-vanishing, there exists a sequence $\{y_k\} \subset \mathbb{R}^N$ and a nontrivial element $v \in H^1(\mathbb{R}^N)$ such that up to a subsequence $v_k(\cdot + y_k) \rightharpoonup v$ in $H^1(\mathbb{R}^N)$ and $v_k( \cdot + y_k) \to v$ almost everywhere on $\mathbb{R}^N$. Denote $m' := \|v\|^2_{L^2(\mathbb{R}^N)} \in (0, m^*]$ and $w_k := v_k(\cdot + y_k) - v$. Clearly, $\lim_{k \to \infty} \|w_k\|^2_{L^2(\mathbb{R}^N)} = m^* -m'$ and
  \begin{equation*}
    \lim_{k \to \infty} \|\nabla v_k\|^2_{L^2(\mathbb{R}^N)} = \|\nabla v\|^2_{L^2(\mathbb{R}^N)} + \lim_{k \to \infty} \|\nabla w_k\|^2_{L^2(\mathbb{R}^N)}.
  \end{equation*}
Noting that $E_{m_k} \to E_{m^*} = 0$ via Theorem \ref{th-1}, we deduce from \eqref{eq:v_k} and the splitting result \cite[Lemma 2.6]{JL20} that
  \begin{equation}\label{eq:Iv_k}
    0 = \lim_{k \to \infty} I(v_k) = \lim_{k \to \infty} I(v + w_k) = I(v) + \lim_{k \to \infty}I(w_k).
  \end{equation}
Since $E_s = 0$ for any $s \in (0, m^*]$ by  Theorem \ref{th-1} $(i)$,  one may see that $I(u) \geq 0$ for any $u \in H^1(\mathbb{R}^N)$ satisfying $\|u\|^2_{L^2(\mathbb{R}^N)} \leq m^*$ and therefore $I(v) \geq 0$ and $\lim_{k \to \infty} I(w_k) \geq 0$. In view also of \eqref{eq:Iv_k}, we obtain
  \begin{equation*}
    I(v) = 0 \qquad \text{and} \qquad \lim_{k \to \infty} I(w_k) = 0.
  \end{equation*}
Now, from the fact that by Theorem \ref{th-1} $(iii)$  the global infimum $E_s = 0$ is not  achieved  for any $s \in (0, m^*)$, it follows that $m' = m^*$ and thus $v \in S_{m^*}$ is a minimizer of $E_{m^*} = 0$. As a consequence, $v \in S_{m^*}$ is an energy ground state solution  of \eqref{problem} when $m = m^*$ and with a Lagrange multiplier $\mu \in \mathbb{R}$. Using the free Pohozaev identity \eqref{eq:pohozaev-libre} it is easy to see that
  \begin{equation*}
    0 = I(v) = I(v) - Q(v) = \frac{1}{N} \|\nabla v \|^2_{L^2(\mathbb{R}^N)} - \frac{1}{2}\mu m^*
  \end{equation*}
and hence $\mu > 0$. The lemma is proved.   \hfill $\square$

Our next lemma is a deformation result which proves convenient to look for normalized radial solutions at positive energy levels via minimax arguments. In particular, it will enable us to apply the genus theory and a suitable version of the Symmetric Mountain Pass Theorem due to \cite{AR73}  to prove Theorem \ref{theorem:multiplicity}. To present the lemma,  we recall  the so-called Palais-Smale-Pohozaev condition at a level $c \in \mathbb{R}$, shortly the $(PSP)_c$ condition, which was introduced in the papers \cite{HT19,IT19}.
\begin{definition}
  For given $c \in \mathbb{R}$, one says that the constrained functional $I_{|S_m \cap H^1_r(\mathbb{R}^N)}$ satisfies the $(PSP)_c$ condition if any sequence $\{u_n\} \subset S_m \cap H^1_r(\mathbb{R}^N)$ for which
    \begin{equation}\label{eq:PSP_c}
      I(u_n) \to c, \qquad  \big( I_{|S_m \cap H^1_r(\mathbb{R}^N)} \big)'(u_n) \to 0  \qquad \text{and} \qquad  P(u_n) \to 0,
    \end{equation}
  has a strongly convergent subsequence in $H^1_r(\mathbb{R}^N)$. For convenience, any sequence $\{u_n\}$ that satisfies \eqref{eq:PSP_c} is hereafter called a $(PSP)_c$ sequence of $I_{|S_m \cap H^1_r(\mathbb{R}^N)}$.
\end{definition}

We also need some notations. For given $u \in H^1(\mathbb{R}^N)$ and any $\theta \in \mathbb{R}$, we define the scaling function
  \begin{equation*}
    (\theta \diamond u)(x) := e^{N\theta/2} u (e^\theta x) \qquad \text{for almost everywhere}~x \in \mathbb{R}^N.
  \end{equation*}
Note that the element $\theta \diamond u \in H^1(\mathbb{R}^N)$ preserves the $L^2$ norm when $\theta \in \mathbb{R}$ varies and that the map sending $(\theta, u) \in \mathbb{R} \times H^1(\mathbb{R}^N)$ to $\theta \diamond u \in H^1(\mathbb{R}^N)$ is continuous. For given $c \in \mathbb{R}$, let
  \begin{equation*}
    I^c:=\{u \in S_m \cap H^1_r(\mathbb{R}^N)~|~I(u) \leq c\}
  \end{equation*}
and denote by  $K^c$ the set of critical points of $I_{|S_m\cap H^1_r(\mathbb{R}^N)}$  at level $c \in \mathbb{R}$.
\begin{lemma}\label{lemma:deformation}
  Assume that $N \geq 2$ and $f \in C(\mathbb{R}, \mathbb{R})$ satisfies $(\hyperlink{f1}{f1}) - (\hyperlink{f2}{f2})$. If the constrained functional $I_{|S_m \cap H^1_r(\mathbb{R}^N)}$ satisfies the $(PSP)_c$ condition at some level $c \in \mathbb{R}$, then for any neighborhood $\mathcal{O} \subset S_m \cap H^1_r(\mathbb{R}^N)$ of $K^c$ ($\mathcal{O} = \emptyset$ if $K^c = \emptyset$) and any $\overline{\varepsilon} > 0$, there exists $\varepsilon \in (0, \overline{\varepsilon})$ and $\eta \in C([0,1] \times (S_m \cap H^1_r(\mathbb{R}^N)), S_m \cap H^1_r(\mathbb{R}^N))$ such that the following properties hold.
    \begin{itemize}
      \item[$(i)$] $\eta(0,u) = u$ for any $u \in S_m \cap H^1_r(\mathbb{R}^N)$.
      \item[$(ii)$] $\eta(t,u) = u$ for any $t \in [0,1]$ if $u \in I^{c- \overline{\varepsilon}}$.
      \item[$(iii)$] $t \mapsto I(\eta(t,u))$ is nonincreasing for any $u \in S_m \cap H^1_r(\mathbb{R}^N)$.
      \item[$(iv)$] $\eta(1, I^{c + \varepsilon} \setminus \mathcal{O}) \subset I^{c - \varepsilon}$ and $\eta(1, I^{c + \varepsilon}) \subset I^{c - \varepsilon} \cup \mathcal{O}$.
      \item[$(v)$] $\eta(t, -u) = - \eta(t,u)$ for any $(t,u) \in [0,1] \times (S_m \cap H^1_r(\mathbb{R}^N))$ when $f$ is odd.
    \end{itemize}
\end{lemma}
\proof We apply directly the abstract deformation result that was developed recently by Ikoma and Tanaka \cite[Proposition 4.5]{IT19}. Indeed, let
  \begin{equation*}
    E = H^1_r(\mathbb{R}^N), \qquad \Phi_\theta u = \theta \diamond u \qquad \text{and} \qquad S = S_m \cap H^1_r(\mathbb{R}^N).
  \end{equation*}
It is not difficult to check that the assumption $(\Phi, S, I)$ made in \cite[Page 637]{IT19} is satisfied and hence the lemma follows.  \hfill $\square$

In order to make use of the above deformation lemma in our later proofs, the $(PSP)_c$ condition must be verified a priori at some relevant level $c \in \mathbb{R}$.
\begin{lemma}\label{lemma:PSP_c}
  Assume that $N \geq 2$ and $f \in C(\mathbb{R}, \mathbb{R})$ satisfies $(\hyperlink{f1}{f1})$, $(\hyperlink{f2}{f2})$ and $(\hyperlink{f5}{f5})$. Then the constrained functional $I_{|S_m \cap H^1_r(\mathbb{R}^N)}$ satisfies the $(PSP)_c$ condition for any $c \neq 0$.
\end{lemma}
\proof Let $\{u_n\} \subset S_m \cap H^1_r(\mathbb{R}^N)$ be an arbitrary $(PSP)_c$ sequence with $c \neq 0$. Since $I$ is coercive on $S_m \cap H^1_r(\mathbb{R}^N)$ by Lemma \ref{lemma:I} $(ii)$, the sequence $\{u_n\}$ is bounded and then there exists $u \in H^1_r(\mathbb{R}^N)$ such that up to a subsequence $u_n \rightharpoonup u$ in $H^1_r(\mathbb{R}^N)$, $u_n \to u$ in $L^p(\mathbb{R}^N)$ for any $p \in (2, 2^*)$, and $u_n \to u$ almost everywhere on $\mathbb{R}^N$.  Moreover, by  \cite[Lemma 3]{Be83-2},
  \begin{equation}\label{eq:u_n}
    - \Delta u_n + \mu_n u_n - f(u_n) \to 0 \qquad \text{in}~\big(H^1_r(\mathbb{R}^N)\big)^{-1},
  \end{equation}
where
  \begin{equation}\label{eq:mu_n}
    \mu_n := \frac{1}{m} \left(\int_{\mathbb{R}^N} f(u_n)u_n dx - \int_{\mathbb{R}^N}|\nabla u_n|^2 dx \right).
  \end{equation}
Since up to a subsequence $\lim_{n \to \infty}\int_{\mathbb{R}^N}|\nabla u_n|^2 dx$ and  $\lim_{n \to \infty} \int_{\mathbb{R}^N} f(u_n)u_n dx$ exist, one may assume that $\mu_n \to \mu$ for some $\mu \in \mathbb{R}$ and thus
  \begin{equation}\label{eq:u}
    - \Delta u + \mu u = f(u) \qquad \text{in}~\big(H^1_r(\mathbb{R}^N)\big)^{-1}.
  \end{equation}
We now claim that
  \begin{equation}\label{eq:nabla-un}
    \lim_{n \to \infty} \|\nabla u_n\|_{L^2(\mathbb{R}^N)} >0.
  \end{equation}
If not, then $\lim_{n \to \infty} \int_{\mathbb{R}^N} F(u_n) dx = 0$ by Lemma \ref{lemma:I} $(i)$ and we obtain thus a contradiction
  \begin{equation*}
    0 \neq c = \lim_{n \to \infty} I(u_n) = 0.
  \end{equation*}
In order to prove that $\mu > 0$, for the number $\vartheta \in (2, 2^*)$ given in $(\hyperlink{f5}{f5})$, we set
  \begin{equation*}
    \beta := \frac{2N - (N - 2)\vartheta}{N(\vartheta - 2)} > 0.
  \end{equation*}
From \eqref{eq:mu_n}, the fact that $P(u_n) \to 0$ and $(\hyperlink{f5}{f5})$, it follows
  \begin{equation*}
    \begin{split}
      m \mu
        & = m \mu_n + o_n(1) \\
        & = \beta \int_{\mathbb{R}^N}|\nabla u_n|^2dx - (1 + \beta) \int_{\mathbb{R}^N}|\nabla u_n|^2dx + \int_{\mathbb{R}^N} f(u_n)u_n dx + o_n(1) \\
        & = \beta \int_{\mathbb{R}^N}|\nabla u_n|^2dx - (1 + \beta) P(u_n) + \frac{2}{\vartheta - 2}\int_{\mathbb{R}^N}\big[\vartheta F(u_n) - f(u_n)u_n \big]dx + o_n(1) \\
        & \geq \beta \int_{\mathbb{R}^N}|\nabla u_n|^2 dx + o_n(1).
    \end{split}
  \end{equation*}
In view of \eqref{eq:nabla-un}, we deduce that $\mu > 0$. To conclude the proof, we denote
  \begin{equation*}
    f_1(t) :=
               \left\{
                  \begin{aligned}
                    & \max \{0, f(t)\}, \qquad & & \text{for} ~ t \geq 0,\\
                    & \min \{0, f(t)\}, \qquad & & \text{for} ~ t < 0,
                  \end{aligned}
               \right.
  \end{equation*}
and $f_2 := f_1 - f$. Since
  \begin{equation*}
    \lim_{t \to 0}\frac{f_1(t)}{t} = 0 = \lim_{t \to \infty} \frac{f_1(t)}{|t|^{1 + \frac{4}{N}}}
  \end{equation*}
and $v_n := u_n - u \to 0$ in $L^{2 + 4/N}(\mathbb{R}^N)$, we have $\int_{\mathbb{R}^N}f_1(u_n)v_ndx \to 0$. Using that $u_n \rightharpoonup u$ in $H^1_r(\mathbb{R}^N)$, it is standard to show also that $\int_{\mathbb{R}^N} \big[f_1(u_n) - f_1(u)\big]u dx \to 0$. Therefore,
  \begin{equation*}
    \lim_{n \to \infty}\int_{\mathbb{R}^N} f_1(u_n)u_ndx = \int_{\mathbb{R}^N} f_1(u)udx.
  \end{equation*}
Since Fatou's lemma implies
  \begin{equation*}
    \int_{\mathbb{R}^N} f_2(u)u dx \leq \lim_{n \to \infty} \int_{\mathbb{R}^N}f_2(u_n)u_n dx
  \end{equation*}
it can be seen that
  \begin{equation}\label{eq:add-1}
    \int_{\mathbb{R}^N} f(u)udx \geq \lim_{n \to \infty}\int_{\mathbb{R}^N} f(u_n)u_ndx.
  \end{equation}
Thus, by \eqref{eq:u_n}, \eqref{eq:u}, \eqref{eq:add-1} and the fact that $\mu_n \to \mu > 0$, we obtain
  \begin{equation*}
    \begin{split}
      \int_{\mathbb{R}^N}|\nabla u|^2dx + \mu \int_{\mathbb{R}^N}|u|^2dx
        & = \int_{\mathbb{R}^N}f(u)udx  \\
        & \geq \lim_{n \to \infty} \int_{\mathbb{R}^N}f(u_n)u_ndx = \lim_{n \to \infty} \int_{\mathbb{R}^N}|\nabla u_n|^2dx + \mu m\\
        & \geq \int_{\mathbb{R}^N}|\nabla u|^2dx + \mu \int_{\mathbb{R}^N}|u|^2dx.
    \end{split}
  \end{equation*}
Clearly, $\lim_{n \to \infty} \int_{\mathbb{R}^N}|\nabla u_n|^2dx = \int_{\mathbb{R}^N}|\nabla u|^2dx$, $\int_{\mathbb{R}^N}|u|^2dx = m$ and so $u_n \to u$ in $H^1_r(\mathbb{R}^N)$.  \hfill $\square$

\begin{remark}\label{rmk:PSP_c}
  The above compactness lemma is valid for any $c \neq 0$, but we shall only use it at suspected positive critical levels of saddle type  in our later arguments. Lemma \ref{lemma:PSP_c} should be compared with \cite[Lemma 5.2]{JL19} which shows that without assuming $(\hyperlink{f5}{f5})$ the constrained functional $I_{|S_m \cap H^1_r(\mathbb{R}^N)}$ satisfies  the Palais-Smale condition at any level $c < 0$, a compactness condition stronger than the $(PSP)_c$ condition.
\end{remark}

%%%%%%%%%%%%%%%%%%%%%%%%%%%%%%%%%%%%%%%%%%%%%%%%%%%%%%%%%%%%%%%%%%%%%%%%%%%%%%%%%%%%%%%%%%%%%%%%%%%%%%%%%%%%%%%%%%%%%%%%%%%%%%%%%%%%%%
%%%%%%%%%%%%%%%%%%%%%%%%%%%%%%%%%%%%%%%%%%%%%%%%%%%%%%%%%%%%%%%%%%%%%%%%%%%%%%%%%%%%%%%%%%%%%%%%%%%%%%%%%%%%%%%%%%%%%%%%%%%%%%%%%%%%%%

\section{Local minimizers}\label{sect:local}

This section is devoted to the proof of Theorem \ref{theorem:local}. For convenience, we recall that
    \begin{equation*}
      \overline{E}_m := \inf_{u \in S^\rho_m} I(u) \geq 0 \qquad \text{for any}~m \in (0, m^*],
    \end{equation*}
where
  \begin{equation*}
    S^\rho_m := \big\{ u \in S_m~|~\|\nabla u\|^2_{L^2(\mathbb{R}^N)} > \rho(m^*) \big\}
  \end{equation*}
and $\rho(m^*) > 0$ is the number given in Lemma \ref{lemma:I} $(iii)$. We shall discuss some properties of the local infimum $\overline{E}_m$ in Subsection \ref{subsect:properties} and then complete the proof of Theorem \ref{theorem:local} in Subsection \ref{subsect:existence}.

%%%%%%%%%%%%%%%%%%%%%%%%%%%%%%%%%%%%%%%%%%%%%%%%%%%%%%%%%%%%%%%%%%%%%%%%%%%%%%%%%%%%%%%%%%%%%%%%%%%%%%%%%%%%%%%%%%%%%%%%%%%%%%%%%%%%%%

\subsection{Some properties of the local infimum}\label{subsect:properties}

In this subsection we put in light some properties of the local infimum $\overline{E}_m$ when the mass $m \in (0, m^*]$ varies. On one hand, this study has its own interest; in particular we manage to show in Lemma \ref{lemma:nonincrease} that $\overline{E}_m$ is nonincreasing, a feature shared by corresponding values in many different $L^2$ constraint problems. On the other hand, as one will see in Lemma \ref{lemma:compactness} of the next subsection, it enables us to develop a robust compactness argument to obtain minimizers of $\overline{E}_m$, which does not require the nonlinearity $f$ to be odd and is valid for any minimizing sequence of the local infimum $\overline{E}_m$. As a first result  we have
\begin{lemma}\label{lemma:continuity}
  Assume that $N\geq1$ and $f\in C(\mathbb{R},\mathbb{R})$ satisfies $(\hyperlink{f1}{f1})-(\hyperlink{f4}{f4})$. Then the function $m \mapsto \overline{E}_m$ is continuous at each $m \in (0, m^*]$.
\end{lemma}
\proof It is sufficient to show that for given $m \in (0, m^*]$ and any sequence $\{m_k\} \subset (0, m^*)$ such that $m_k \to m$ as $k \to \infty$ one has $\lim_{k \to \infty} \overline{E}_{m_k} = \overline{E}_m$. We first prove
  \begin{equation}\label{eq:upper}
    \limsup_{k \to \infty} \overline{E}_{m_k} \leq \overline{E}_m.
  \end{equation}
For any $u \in S^\rho_m$ and each $k \in \mathbb{N}^+$, set $u_k:= \sqrt{m_k/m} \cdot u \in S_{m_k}$. Since $u_k \to u$ in $H^1(\mathbb{R}^N)$, it is clear that $u_k \in S^\rho_{m_k}$ for any $k$ large enough and $\lim_{k \to \infty} I(u_k) = I(u)$. Thus
  \begin{equation*}
    \limsup_{k \to \infty} \overline{E}_{m_k} \leq \limsup_{k \to \infty} I(u_k) = I(u).
  \end{equation*}
By the arbitrariness of $u \in S^\rho_m$, we conclude that \eqref{eq:upper} holds. To complete the proof, it remains to show
  \begin{equation}\label{eq:lower}
    \liminf_{k \to \infty} \overline{E}_{m_k} \geq \overline{E}_m.
  \end{equation}
For each $k \in \mathbb{N}^+$, there exists $v_k \in S^\rho_{m_k}$ such that
  \begin{equation}\label{eq:lower-1}
    I(v_k) \leq \overline{E}_{m_k} + \frac{1}{k}.
  \end{equation}
Setting $t_k := \sqrt{m/m_k}$, we have $\tilde{v}_k := t^{(2 - N)/2}_k v_k(\cdot/t_k) \in S^\rho_m$ and thus
  \begin{equation*}
    \begin{split}
      \overline{E}_m \leq I(\tilde{v}_k)
        & \leq I(v_k) + \big|I(\tilde{v}_k) - I(v_k)\big|\\
        & \leq \overline{E}_{m_k} + \frac{1}{k} + \big|I(\tilde{v}_k) - I(v_k)\big| \\
        & \leq \overline{E}_{m_k} + \frac{1}{k} +  t^N_k \int_{\mathbb{R}^N} \left|F\Big(t^{\frac{2 - N}{2}}_k v_k\Big)- F(v_k)\right| dx + |t^N_k - 1| \int_{\mathbb{R}^N}|F(v_k)|dx.
    \end{split}
  \end{equation*}
Since $t_k \to 1$ and $f \in C(\mathbb{R}, \mathbb{R})$ satisfies $(\hyperlink{f1}{f1})$ and $(\hyperlink{f2}{f2})$, the proof of \eqref{eq:lower} can be reduced to show that $\{v_k\}$ is bounded in $H^1(\mathbb{R}^N)$. To justify the boundedness, by \eqref{eq:lower-1} and \eqref{eq:upper}, we have $\limsup_{k \to \infty} I(v_k) \leq \overline{E}_m$. Noting that $v_k \in S_{m_k}$ and $m_k \to m$, it follows from Lemma \ref{lemma:I} $(ii)$ that $\{v_k\}$ is bounded in $H^1(\mathbb{R}^N)$. \hfill $\square$

\begin{lemma}\label{lemma:nonincrease}
  Assume that $N\geq1$ and $f\in C(\mathbb{R},\mathbb{R})$ satisfies $(\hyperlink{f1}{f1})-(\hyperlink{f4}{f4})$. Then the function $m \mapsto \overline{E}_m$ is nonincreasing on $(0, m^*]$.
\end{lemma}
\proof It is equivalent to show that for any $m, m' \in (0, m^*]$ satisfying $m > m'$ and for an arbitrary $\varepsilon > 0$ one has
  \begin{equation}\label{eq:nonincrease_1}
    \overline{E}_m \leq \overline{E}_{m'} + \varepsilon.
  \end{equation}
By the definition of $\overline{E}_{m'}$, there exists $u \in S^\rho_{m'}$ such that
  \begin{equation}\label{eq:nonincrease_2}
    I(u) \leq \overline{E}_{m'} + \frac{\varepsilon}{2}.
  \end{equation}
Let $\chi\in C^\infty_c(\mathbb{R}^N)$ be a radial cut-off function such that $\chi(x) = 1$ for $|x| \leq 1$ and $\chi(x) = 0$ when $|x| \geq 2$. For any $\delta > 0$, we set $u_\delta(x) := u(x) \chi(\delta x)$. Since $u_\delta \to u$ in $H^1(\mathbb{R}^N)$ as $\delta \to 0^+$, one can fix a small enough constant $\delta > 0$ such that $\| \nabla u_\delta \|^2_{L^2(\mathbb{R}^N)} > \rho(m^*)$ and
  \begin{equation}\label{eq:nonincrease_3}
    I(u_\delta) \leq I(u) + \frac{\varepsilon}{4}.
  \end{equation}
Then take $v \in C^\infty_c(\mathbb{R}^N) \setminus \{0\}$ such that $\text{supp}(v) \subset B(0,1+4/\delta) \setminus B(0,4/\delta)$ and set
\begin{equation*}
    \tilde{v} := \frac{\sqrt{m - \|u_\delta\|^2_{L^2(\mathbb{R}^N)}}}{\|v\|_{L^2(\mathbb{R}^N)}} v.
  \end{equation*}
For any $\lambda \leq 0$, we define $w_\lambda := u_\delta + \lambda \diamond \tilde{v}$. Since $u_\delta$ and  $\lambda \diamond \tilde{v}$ have disjoint supports, it is clear that $w_\lambda \in S^\rho_m$.  Noting that $\|\nabla (\lambda \diamond \tilde{v})\|_{L^2(\mathbb{R}^N)} \to 0$ as $\lambda \to - \infty$, it follows from Lemma \ref{lemma:I} $(i)$ that
  \begin{equation}\label{eq:nonincrease_4}
    I(\lambda_0 \diamond \tilde{v}) \leq \frac{\varepsilon}{4} \qquad \text{for some}~\lambda_0 < 0.
  \end{equation}
Now, by the definition of $\overline{E}_m$, \eqref{eq:nonincrease_4}, \eqref{eq:nonincrease_3} and \eqref{eq:nonincrease_2}, we obtain
  \begin{equation*}
    \overline{E}_m \leq I(w_{\lambda_0}) = I(u_\delta) + I(\lambda_0 \diamond \tilde{v}) \leq I(u) + \frac{\varepsilon}{2} \leq \overline{E}_{m'} + \varepsilon,
  \end{equation*}
that is \eqref{eq:nonincrease_1}.  \hfill $\square$
\medskip

When the local infimum $\overline{E}_{m'}$ is  achieved  by some $u \in S^\rho_{m'}$ with $m' \in (0, m^*]$, knowing the sign of the corresponding Lagrange multiplier  provides accurate information on the monotonicity of $\overline{E}_m$ at the mass $m = m'$.
\begin{lemma}\label{lemma:monotonicity}
  Assume that $N \geq 1$ and $f \in C(\mathbb{R}, \mathbb{R})$ satisfies $(\hyperlink{f1}{f1})-(\hyperlink{f4}{f4})$. Suppose that for some $m' \in (0, m^*]$ there exists a couple $(u, \mu) \in S^\rho_{m'} \times \mathbb{R}$  such that
    \begin{equation*}
      - \Delta u + \mu u = f(u)
    \end{equation*}
  and $I(u) = \overline{E}_{m'}$. Then $\overline{E}_m < \overline{E}_{m'}$ for any $m < m'$ close enough to $m'$ if $\mu < 0$ and for each $m \in (m', m^*]$ near enough to $m'$ if $\mu > 0$.
\end{lemma}
\proof Let $u \in S^\rho_{m'}$ and $\mu \in \mathbb{R}$ be as above. For any $t > 0$, we set $u_t := tu \in S_{m' t^2}$ and
  \begin{equation*}
    \alpha(t) := I(u_t) = \frac{1}{2} t^2 \int_{\mathbb{R}^N} |\nabla u|^2 dx - \int_{\mathbb{R}^N} F(tu) dx.
  \end{equation*}
When $\mu < 0$, by the facts that $u_t \to u$ strongly in $H^1(\mathbb{R}^N)$ as $t \to 1$ and that
  \begin{equation*}
    I'(u)u = -\mu \|u\|^2_{L^2(\mathbb{R}^N)} = - \mu m' > 0,
  \end{equation*}
one may fix a small enough constant $\delta > 0$ such that for any $t \in [1 - \delta, 1)$,
  \begin{equation*}
    u_t \in S^\rho_{m' t^2} \qquad \text{and} \qquad \frac{d}{dt} \alpha(t) = t^{-1} I'(u_t)u_t > 0.
  \end{equation*}
Then, from the mean value theorem, it follows that
  \begin{equation*}
    \overline{E}_{m' t^2} \leq \alpha(t) = \alpha(1) + (t-1) \cdot \frac{d}{dt} \alpha(\theta) < \alpha(1),
  \end{equation*}
where $1 - \delta \leq t < \theta < 1$. For any $m < m'$ close enough to $m'$, we have
  \begin{equation*}
    t := \sqrt{m/m'} \in [1 - \delta, 1)
  \end{equation*}
and thus $\overline{E}_m < \alpha(1) = I(u) = \overline{E}_{m'}$. The case of $\mu > 0$ can be proved similarly.  \hfill $\square$
\medskip

From Lemmas \ref{lemma:nonincrease} and \ref{lemma:monotonicity}, we directly obtain
\begin{lemma}\label{lemma:decrease}
  Assume $N\geq1$ and $f\in C(\mathbb{R},\mathbb{R})$ satisfies $(\hyperlink{f1}{f1})-(\hyperlink{f4}{f4})$. If for some $m' \in (0, m^*]$ there exists a couple $(u, \mu) \in S^\rho_{m'} \times \mathbb{R}$  such that
    \begin{equation*}
      - \Delta u + \mu u = f(u)
    \end{equation*}
  and $I(u) = \overline{E}_{m'}$, then $\mu \geq 0$. If in addition $\mu > 0$, then $\overline{E}_m < \overline{E}_{m'}$ for any $m \in (m', m^*]$.
\end{lemma}

%%%%%%%%%%%%%%%%%%%%%%%%%%%%%%%%%%%%%%%%%%%%%%%%%%%%%%%%%%%%%%%%%%%%%%%%%%%%%%%%%%%%%%%%%%%%%%%%%%%%%%%%%%%%%%%%%%%%%%%%%%%%%%%%%%%%%%

\subsection{Existence of local minimizers}\label{subsect:existence}

In this subsection we  show the existence of local minimizers for a suitable range of the mass and then finish the proof of Theorem \ref{theorem:local}. We establish at first the following geometrical result about the local infimum $\overline{E}_m$.
\begin{lemma}\label{lemma:local}
  Assume that $N \geq 1$ and $f \in C(\mathbb{R},\mathbb{R})$ satisfies $(\hyperlink{f1}{f1})-(\hyperlink{f4}{f4})$. Then there exists $m^{**} \in (0, m^*)$ such that for any $m \in (m^{**}, m^*]$
    \begin{equation*}
      \overline{E}_{m} := \inf_{u \in S^\rho_m} I(u) < \Lambda := \min\left\{\frac{1}{4}, \frac{1}{N}\right\} \rho(m^*) \leq \inf_{u \in \mathcal{A}_m} I(u),
    \end{equation*}
  where $\mathcal{A}_m := \{u \in S_m~|~\rho(m^*) < \|\nabla u\|^2_{L^2(\mathbb{R}^N)} \leq 4 \rho(m^*)\}$. In particular, $\overline{E}_{m^*} = 0$.
\end{lemma}
\proof By Lemma \ref{lemma:m*}, we can fix a minimizer $v \in S_{m^*}$ of $E_{m^*} = 0$. From Lemma \ref{lemma:I} $(iii)$, it follows that
  \begin{equation*}
    \|\nabla v\|^2_{L^2(\mathbb{R}^N)} > 4 \rho(m^*) \qquad \text{and} \qquad \inf_{u \in \mathcal{A}_m} I(u) \geq \frac{1}{4}\rho(m^*).
  \end{equation*}
Then, by continuity there exists $\tau \in (0, 1)$ such that for any $t \in (\tau, 1]$
  \begin{equation*}
    \|\nabla (t v)\|^2_{L^2(\mathbb{R}^N)} > 4\rho(m^*) \qquad \text{and} \qquad I(t v) < \min\left\{\frac{1}{4}, \frac{1}{N}\right\} \rho(m^*).
  \end{equation*}
Clearly, the proof is complete with the choice of $m^{**} := \tau^2 m^*$. \hfill $\square$
\begin{remark}
  Being the minimum of $\rho(m^*)/4$ and $\rho(m^*)/N$, the upper bound $\Lambda$ on $\overline{E}_m$ has a threefold role in our later proofs. First, it allows us to rule out the vanishing case for any minimizing sequence $\{u_n\} \subset S^\rho_m$ of $\overline{E}_m$. Second, Lemma \ref{lemma:I} $(iii)$ gives that $\|\nabla u\|^2_{L^2(\mathbb{R}^N)} > 4 \rho(m^*)$ for any $u \in S^\rho_m$ with $I(u) < \rho(m^*)/4$, and hence an arbitrary minimizing sequence $\{u_n\} \subset S^\rho_m$ of $\overline{E}_m$ satisfies
    \begin{equation*}
      \|\nabla u_n\|^2_{L^2(\mathbb{R}^N)} > 4 \rho(m^*) \qquad\text{for any}~n~\text{large enough}.
    \end{equation*}
  As one will see in Lemma \ref{lemma:compactness} below,  such an extra information is crucial to overcome certain difficulties caused by the local constraint $\|\nabla u\|^2_{L^2(\mathbb{R}^N)} > \rho(m^*)$ on $S_m$. Last but not least, the upper estimate $\overline{E}_m < \rho(m^*)/N$ permits to show that the Lagrange multiplier  in \eqref{eq:solution} below  is strictly positive, and this proves important in our compactness argument when we look for minimizers of the local infimum $\overline{E}_m$.
\end{remark}

In order to address the compactness issue when looking for minimizers of $\overline{E}_m$, we develop below a compactness argument which does not require the nonlinearity $f$ to be odd and is valid for any minimizing sequence of $\overline{E}_m$.
\begin{lemma}\label{lemma:compactness}
  Assume that $N \geq 1$ and $f \in C(\mathbb{R}, \mathbb{R})$ satisfies $(\hyperlink{f1}{f1})-(\hyperlink{f4}{f4})$. Let $m \in (m^{**}, m^*]$ and $\{v_n\} \subset S^\rho_m$ be an arbitrary minimizing sequence of $\overline{E}_m$. Then up to a subsequence there exists a sequence $\{y_n\} \subset \mathbb{R}^N$ and a nontrivial element $v \in S^\rho_m$ such that $v_n( \cdot + y_n)$ converges strongly to $v$ in $H^1(\mathbb{R}^N)$. In particular, $v \in S^\rho_m$ is a minimizer of the local infimum $\overline{E}_m$ with
    \begin{equation}\label{eq:minimum}
      I(v) =
        \left\{
          \begin{aligned}
          & ~ \overline{E}_m > 0, \qquad \text{for}~m \in (m^{**}, m^*),\\
          & ~ \overline{E}_m = 0, \qquad \text{for}~m = m^*,
        \end{aligned}
      \right.
    \end{equation}
  and there exists a positive Lagrange multiplier $\mu > 0$ such that $- \Delta v + \mu v = f(v)$.
\end{lemma}
\proof Since $I_{|S_m}$ is coercive by Lemma \ref{lemma:I} $(ii)$,  the minimizing sequence $\{v_n\} \subset S^\rho_m$ is bounded in $H^1(\mathbb{R}^N)$ and thus one may assume that up to a subsequence $\lim_{n \to \infty} \int_{\mathbb{R}^N} |\nabla v_n|^2 dx$ and $\lim_{n \to \infty} \int_{\mathbb{R}^N}F(v_n) dx$ exist.  We claim that $\{v_n\}$ is non-vanishing, that is
  \begin{equation*}
    \sigma := \limsup_{n \to \infty} \left(\sup_{y \in \mathbb{R}^N} \int_{B(y, 1)}|v_n|^2 dx \right) > 0.
  \end{equation*}
Indeed, if $\{v_n\}$ were vanishing, namely $\sigma = 0$, then \cite[Lemma I.1]{Lions84-2} would imply that $\|v_n\|_{L^{2 + 4/N}(\mathbb{R}^N)} \to 0$. In view of Lemma \ref{lemma:I} $(i)$, we have
  \begin{equation*}
    \overline{E}_m = \lim_{n \to \infty} I(v_n) \geq \frac{1}{2}\lim_{n \to \infty}\|\nabla v_n\|^2_{L^2(\mathbb{R}^N)} \geq \frac{1}{2}\rho(m^*).
  \end{equation*}
This contradicts Lemma \ref{lemma:local} and thus the claim is proved.

Since $\{v_n\}$ is non-vanishing, there exists a sequence $\{y_n\} \subset \mathbb{R}^N$ and a nontrivial element $v \in H^1(\mathbb{R}^N)$ such that up to a subsequence $v_n(\cdot + y_n) \rightharpoonup v$ in $H^1(\mathbb{R}^N)$ and $v_n(\cdot + y_n) \to v$ almost everywhere on $\mathbb{R}^N$. Denote $m' := \|v\|^2_{L^2(\mathbb{R}^N)} \in (0, m]$ and $w_n := v_n( \cdot + y_n) - v$. It is clear that $\lim_{n \to \infty} \|w_n\|^2_{L^2(\mathbb{R}^N)} = m - m'$,
  \begin{equation}\label{eq:BL-nabla}
    \lim_{n \to \infty} \|\nabla v_n\|^2_{L^2(\mathbb{R}^N)} = \|\nabla v \|^2_{L^2(\mathbb{R}^N)} + \lim_{n \to \infty} \|\nabla w_n\|^2_{L^2(\mathbb{R}^N)}
  \end{equation}
and
  \begin{equation}\label{eq:BL-I}
    \overline{E}_m = \lim_{n \to \infty} I(v_n) = \lim_{n \to \infty} I(v + w_n)= I(v) + \lim_{n \to \infty} I(w_n).
  \end{equation}
One should note that the splitting result \cite[Lemma 2.6]{JL20} was used in \eqref{eq:BL-I}.

The next crucial task is to locate the weak limit $v \in S_{m'}$ on the local constraint $S^\rho_{m'}$, namely to show that
  \begin{equation}\label{eq:nabla-v}
    \|\nabla v\|^2_{L^2(\mathbb{R}^N)} > \rho(m^*).
  \end{equation}
If \eqref{eq:nabla-v} were not true, then $\|\nabla v\|^2_{L^2(\mathbb{R}^N)} \leq \rho(m^*)$. Since $\lim_{n \to \infty} \|\nabla v_n\|^2_{L^2(\mathbb{R}^N)} \geq 4 \rho(m^*)$ by Lemma \ref{lemma:local}, it follows from \eqref{eq:BL-nabla} that
  \begin{equation}\label{eq:nabla-wn}
    \|\nabla w_n\|^2_{L^2(\mathbb{R}^N)} \geq 2 \rho(m^*) \qquad \text{for any} ~n~\text{large enough}.
  \end{equation}
To obtain a contradiction, we distinguish the two cases: {\it compactness} and {\it non-compactness}.

$\bullet$ {\it Compactness: that is $m' = m$.} Then $\|w_n\|_{L^{2 + 4/N}(\mathbb{R}^N)} \to 0$. In view of Lemma \ref{lemma:I} $(i)$ and \eqref{eq:nabla-wn}, it is clear that
  \begin{equation*}
    \lim_{n \to \infty}I(w_n)  \geq \frac{1}{2}\lim_{n \to \infty}\|\nabla w_n\|^2_{L^2(\mathbb{R}^N)} \geq \rho(m^*).
  \end{equation*}
Using \eqref{eq:BL-I} and the fact that $I(v) \geq E_{m'} = 0$ by  Theorem \ref{th-1} $(i)$,  we have
  \begin{equation*}
    \overline{E}_m = I(v) + \lim_{n \to \infty} I(w_n) \geq \rho(m^*),
  \end{equation*}
which contradicts Lemma \ref{lemma:local}.

$\bullet$ {\it Non-compactness: that is $m' < m$.} In this case, $t_n := \|w_n\|^2_{L^2(\mathbb{R}^N)} \to m - m' \in (0, m)$ and thus, by \eqref{eq:nabla-wn}, $w_n \in S^\rho_{t_n}$ for any $n$ large enough. Using the fact that $\overline{E}_m$ is nonincreasing by Lemma \ref{lemma:nonincrease}, we have
  \begin{equation*}
    \lim_{n \to \infty} I(w_n) \geq \limsup \overline{E}_{t_n} \geq \overline{E}_m.
  \end{equation*}
Since $E_{m'} = 0$ is not  achieved by Theorem \ref{th-1} $(iii)$,  it follows that $I(v) > E_{m'} = 0$. In view of \eqref{eq:BL-I}, we have
  \begin{equation*}
    \overline{E}_m = I(v) + \lim_{n \to \infty} I(w_n) > \overline{E}_m,
  \end{equation*}
which is a contradiction.

With the desired location estimate \eqref{eq:nabla-v} at hand, we directly obtain
  \begin{equation*}
    I(v) \geq \overline{E}_{m'}.
  \end{equation*}
Since $E_s = 0$ for any $s \in (0, m^*]$ by  Theorem \ref{th-1} $(i)$,  it is clear that $I(u) \geq 0$ for any $u \in H^1(\mathbb{R}^N)$ satisfying $\|u\|^2_{L^2(\mathbb{R}^N)} \leq m^*$ and thus $\lim_{n \to \infty} I(w_n) \geq 0$. Using \eqref{eq:BL-I} and the fact that $\overline{E}_m$ is nonincreasing by Lemma \ref{lemma:nonincrease}, it is easy to see that
  \begin{equation}\label{eq:Iwn}
    \lim_{n \to \infty}I(w_n) = 0
  \end{equation}
and
  \begin{equation}\label{eq:Iv}
    \overline{E}_m = I(v) = \overline{E}_{m'}.
  \end{equation}
In particular, the weak limit $v \in S^\rho_{m'}$ is a minimizer of $\overline{E}_{m'}$ and hence there exists a Lagrange multiplier $\mu \in \mathbb{R}$ such that
  \begin{equation}\label{eq:solution}
    - \Delta v + \mu v = f(v)
  \end{equation}
	and also $Q(v)=0$, where $Q$ is defined in \eqref{eq:pohozaev-libre}. Then, by Lemma \ref{lemma:local} and \eqref{eq:Iv},
  \begin{equation*}
    \frac{1}{N} \rho(m^*) > \overline{E}_m = I(v) = I(v) - Q(v) = \frac{1}{N} \|\nabla v \|^2_{L^2(\mathbb{R}^N)} - \frac{1}{2} \mu \|v\|^2_{L^2(\mathbb{R}^N)}.
  \end{equation*}
Taking \eqref{eq:nabla-v} into account we conclude that $\mu >0$.

Now, in view of \eqref{eq:Iv}, \eqref{eq:solution} and Lemma \ref{lemma:decrease}, it is clear that $m' = m$ and thus $v \in S^\rho_m$ is a minimizer of $\overline{E}_m$. In particular, \eqref{eq:minimum} follows from the facts that $E_m = 0$ is not  achieved  when $m \in (0, m^*)$ and that $\overline{E}_{m^*} = 0$ by Lemma \ref{lemma:local}. In order to show the strong convergence, noting that $\|w_n\|^2_{L^2(\mathbb{R}^N)} \to m - m' = 0$, we have $\|w_n\|_{L^{2 + 4/N}(\mathbb{R}^N)} \to 0$. In view of \eqref{eq:Iwn} and Lemma \ref{lemma:I} $(i)$,
  \begin{equation*}
    \lim_{n \to \infty}\|\nabla w_n\|^2_{L^2(\mathbb{R}^N)} = 2 \lim_{n \to \infty} \left( I(w_n) + \int_{\mathbb{R}^N} F(w_n) dx \right) \leq 0,
  \end{equation*}
and hence $v_n (\cdot + y_n)$ converges strongly to $v \in S^\rho_m$ in $H^1(\mathbb{R}^N)$.   \hfill $\square$

The next lemma gives information on the sign of minimizers of $\overline{E}_m$ when $m \in (m^{**}, m^*]$. It is valid for all minimizers of $\overline{E}_m$ and does not require the nonlinearity $f$ to be odd.
\begin{lemma}\label{lemma:sign}
  Assume that $N \geq 1$ and $f \in C(\mathbb{R}, \mathbb{R})$ satisfies $(\hyperlink{f1}{f1})-(\hyperlink{f4}{f4})$. Then when $m \in (m^{**}, m^*]$ any minimizer of $\overline{E}_m$ has a constant sign.
\end{lemma}
\proof  For given minimizer $v \in S^\rho_m$ of $\overline{E}_m$, we set
  \begin{equation*}
    v^+ := \max\{0, v\} \qquad \text{and} \qquad v^- := \min\{0, v\},
  \end{equation*}
and suppose by contradiction that $m^+ := \|v^+\|^2_{L^2(\mathbb{R}^N)} \neq 0$ and $m^- := \|v^-\|^2_{L^2(\mathbb{R}^N)} \neq 0$.  Clearly, Lemma \ref{lemma:local} gives
  \begin{equation*}
    4 \rho(m^*) < \|\nabla v\|^2_{L^2(\mathbb{R}^N)} = \|\nabla v^+ \|^2_{L^2(\mathbb{R}^N)} + \|\nabla v^-\|^2_{L^2(\mathbb{R}^N)}.
  \end{equation*}
Without loss of generality, we may assume that $v^+ \in S^\rho_{m^+}$ and thus $I(v^+) \geq \overline{E}_{m^+}$. Since by Theorem \ref{th-1} $(iii)$ the global infimum $E_{m^-} = 0$ is not achieved, $I(v^-) > E_{m^-} = 0$. Using also Lemma \ref{lemma:nonincrease}, we obtain a contradiction
  \begin{equation*}
    \overline{E}_m = I(v) = I(v^+) + I(v^-) > \overline{E}_{m^+} \geq \overline{E}_m.
  \end{equation*}
Hence $v$ has a constant sign on $\mathbb{R}^N$. \hfill $\square$

Having Lemmas \ref{lemma:compactness} and \ref{lemma:sign} at disposal,  we are now able to prove Theorem \ref{theorem:local}.

\medskip
\noindent
{\bf Proof of Theorem \ref{theorem:local}.}  Let $m^{**} \in (0, m^*)$ be the number given in Lemma \ref{lemma:local}. It is clear that most of the conclusions follow from  Lemmas \ref{lemma:continuity}, \ref{lemma:nonincrease}, \ref{lemma:decrease}, \ref{lemma:compactness} and \ref{lemma:sign}.  When \eqref{eq:key_3} also holds, since any solution $u$ of \eqref{problem} satisfies the Pohozaev identity $P(u) = 0$, we conclude from Lemma \ref{lemma:I} $(iii)$ that $u \in S^\rho_m$ and hence the local minimizer $v \in S^\rho_m$ is indeed an energy ground state solution of \eqref{problem}. \hfill $\square$

%%%%%%%%%%%%%%%%%%%%%%%%%%%%%%%%%%%%%%%%%%%%%%%%%%%%%%%%%%%%%%%%%%%%%%%%%%%%%%%%%%%%%%%%%%%%%%%%%%%%%%%%%%%%%%%%%%%%%%%%%%%%%%%%%%%%%%
%%%%%%%%%%%%%%%%%%%%%%%%%%%%%%%%%%%%%%%%%%%%%%%%%%%%%%%%%%%%%%%%%%%%%%%%%%%%%%%%%%%%%%%%%%%%%%%%%%%%%%%%%%%%%%%%%%%%%%%%%%%%%%%%%%%%%%

\section{Mountain pass solutions}\label{sect:MP}

This section aims to prove Theorem \ref{theorem:MP} by using mountain pass arguments. Due to a technical reason, we denote $\varrho(m) := \rho(m^*)$ when $m \in (0, m^*)$ and $\varrho(m) := \rho(m)$ if $m \geq m^*$, where $\rho(m) > 0$ is the number  given in  Lemma \ref{lemma:I} $(iii)$. For any $m > 0$, we introduce the set of continuous paths
  \begin{equation*}
    \Gamma_m :=
      \left\{
        \gamma \in C([0, 1], S_m\cap H^1_r(\mathbb{R}^N)) \left| \begin{array}{c}
                                                               \|\nabla \gamma(0)\|^2_{L^2(\mathbb{R}^N)} < \varrho(m), \\
                                                               \|\nabla \gamma(1)\|^2_{L^2(\mathbb{R}^N)} >  4 \varrho(m), \\
                                                               \max \big\{ I(\gamma(0)), I(\gamma(1)) \big\} < \frac{1}{2} \varrho(m)
                                                            \end{array}
        \right.
      \right\}.
  \end{equation*}
Being the basis for our later discussions, the nonemptiness of $\Gamma_m$ is proved below.
\begin{lemma}\label{lemma:Gamma_m}
  Assume that $N \geq 2$ and $f \in C(\mathbb{R}, \mathbb{R})$ satisfies $(\hyperlink{f1}{f1})-(\hyperlink{f4}{f4})$. Then $\Gamma_m \neq \emptyset$ for any $m > m^{**}$.
\end{lemma}
\proof We first observe  that, if we find $u \in S_m \cap H^1_r(\mathbb{R}^N)$ satisfying
  \begin{equation}\label{eq:test}
    \|\nabla u\|^2_{L^2(\mathbb{R}^N)} >  4 \varrho(m) \qquad \text{and} \qquad I(u) < \frac{1}{2}\varrho(m),
  \end{equation}
then $\Gamma_m \neq \emptyset$. Indeed, for such $u \in S_m \cap H^1_r(\mathbb{R}^N)$,  in view of the fact that $\| \nabla (s \diamond u)\|_{L^2(\mathbb{R}^N)}$ converges to zero as $s \to - \infty$, one may choose $\lambda < 0$  such that
  \begin{equation*}
    \| \nabla (\lambda \diamond u) \|^2_{L^2(\mathbb{R}^N)} < \frac{1}{4}\varrho(m)
  \end{equation*}
and, by Lemma \ref{lemma:I} $(i)$,  $| \int_{\mathbb{R}^N} F(\lambda \diamond u) dx | < \varrho(m)/4$. As a consequence,
  \begin{equation*}
    I(\lambda \diamond u) < \frac{1}{2}\varrho(m).
  \end{equation*}
Setting $\gamma(t) := (\lambda(1 - t)) \diamond u$ for any $t \in [0, 1]$ and combining \eqref{eq:test}, we see that $\gamma \in \Gamma_m$.

To complete this lemma, we only need to look for an element $u \in S_m \cap H^1_r(\mathbb{R}^N)$ that satisfies \eqref{eq:test}. When $m \geq m^*$, by Theorem \ref{th-1} $(ii)$  and Lemma \ref{lemma:m*}, the constrained functional $I_{|S_m}$ admits a global minimizer $v \in S_m$ with $I(v) \leq 0$. In view of  \cite[Theorem 1.4]{JL21},  this minimizer is radially symmetric up to a translation in $\mathbb{R}^N$ and hence we may assume that $v \in S_m \cap H^1_r(\mathbb{R}^N)$. Since $I(v) \leq 0$, it follows from Lemma \ref{lemma:I} $(iii)$ that
  \begin{equation*}
    \|\nabla v\|^2_{L^2(\mathbb{R}^N)} > 4 \varrho(m)
  \end{equation*}
and thus $u := v \in S_m \cap H^1_r(\mathbb{R}^N)$ satisfies \eqref{eq:test}. When $m \in (m^{**}, m^*)$,  since the global minimizer $v \in S_{m^*}$ is radial,  it is clear that $u := \sqrt{m/m^*} \cdot v \in S_m \cap H^1_r(\mathbb{R}^N)$ satisfies \eqref{eq:test}. \hfill $\square$

\begin{lemma}\label{lemma:E_mp}
  Assume that $N \geq 2$ and $f \in C(\mathbb{R}, \mathbb{R})$ satisfies $(\hyperlink{f1}{f1})-(\hyperlink{f4}{f4})$. For any $m > m^{**}$ we define the mountain pass value
    \begin{equation*}
      E_{mp, m} := \inf_{\gamma \in \Gamma_m} \max_{t \in [0, 1]} I(\gamma(t)).
    \end{equation*}
  Then $E_{mp, m} \geq \varrho(m) >0$.
\end{lemma}
\proof Since $\Gamma_m$ is nonempty by Lemma \ref{lemma:Gamma_m}, the mountain pass value $E_{mp, m}$ is well defined. For any given $\gamma \in \Gamma_m$, we have
  \begin{equation*}
    \|\nabla \gamma(0)\|^2_{L^2(\mathbb{R}^N)} < \varrho(m) \qquad \text{and} \qquad \|\nabla \gamma(1)\|^2_{L^2(\mathbb{R}^N)} > 4 \varrho(m),
  \end{equation*}
and then by the intermediate value theorem there exists $\tau \in (0, 1)$ such that
  \begin{equation*}
    \|\nabla \gamma(\tau)\|^2_{L^2(\mathbb{R}^N)} = 4 \varrho(m).
  \end{equation*}
In view of the definition of $\varrho(m)$ and Lemma \ref{lemma:I} $(iii)$, we obtain
  \begin{equation*}
    \max_{t \in [0, 1]} I(\gamma(t)) \geq I(\gamma(\tau)) \geq \varrho(m)
  \end{equation*}
and hence $E_{mp, m} \geq \varrho(m) >0$. \hfill $\square$

Using the deformation result Lemma \ref{lemma:deformation} and the compactness result Lemma \ref{lemma:PSP_c}, we can now obtain  a normalized radial solution at the  mountain pass level.
\begin{lemma}\label{lemma:solution}
  Assume that $N \geq 2$ and $f \in C(\mathbb{R}, \mathbb{R})$ satisfies $(\hyperlink{f1}{f1})-(\hyperlink{f5}{f5})$. Then for each $m > m^{**}$ the constrained functional $I_{|S_m\cap H^1_r(\mathbb{R}^N)}$ admits a radial solution at the mountain pass level $E_{mp, m}$.
\end{lemma}
\proof Let $c := E_{mp, m} \geq \varrho(m) > 0$ and suppose by contradiction that $K^c = \emptyset$. Applying the deformation result Lemma \ref{lemma:deformation} with $\mathcal{O} = \emptyset$ and $\overline{\varepsilon} = \varrho(m)/2 > 0$, there exists $\varepsilon \in (0, \overline{\varepsilon})$ and $\eta \in C([0,1] \times (S_m \cap H^1_r(\mathbb{R}^N)), S_m \cap H^1_r(\mathbb{R}^N))$ such that
  \begin{equation}\label{eq:deformation_MP}
    \eta(1, I^{c + \varepsilon}) \subset I^{c - \varepsilon}.
  \end{equation}
By the definition of $E_{mp, m}$ we can choose $\gamma \in \Gamma_m$ such that
  \begin{equation}\label{eq:gamma}
    \max_{t \in [0, 1]} I(\gamma(t)) \leq c + \varepsilon
  \end{equation}
and consider the new path $\overline{\gamma}(t) := \eta(1, \gamma(t))$ for $t \in [0, 1]$. Since
  \begin{equation*}
    \max \big\{ I(\gamma(0)), I(\gamma(1)) \big\} < \frac{1}{2} \varrho(m) \leq c - \overline{\varepsilon},
  \end{equation*}
it follows from Lemma \ref{lemma:deformation} $(ii)$ that $\overline{\gamma} \in \Gamma_m$. Now, by the definition of $E_{mp, m}$, \eqref{eq:deformation_MP} and \eqref{eq:gamma}, we obtain $c \leq \max_{t \in [0, 1]} I(\overline{\gamma}(t)) \leq c - \varepsilon$, which is a contradiction.   \hfill $\square$

\begin{remark}\label{rmk:positive}
  When the nonlinearity $f$ is odd, we set
    \begin{equation*}
      \widetilde{f}(t) :=
                         \left\{
                           \begin{aligned}
                             & f(t), & \quad \text{for} ~ t \geq 0,\\
                             & ~~ 0, & \quad \text{for} ~ t < 0,
                           \end{aligned}
                         \right.
    \end{equation*}
  denote by $\widetilde{I}$ the associated energy functional, and introduce the new set of continuous paths
    \begin{equation*}
      \widetilde{\Gamma}_m :=
        \left\{
          \gamma \in C([0, 1], S_m\cap H^1_r(\mathbb{R}^N)) \left| \begin{array}{c}
                                                               \|\nabla \gamma(0)\|^2_{L^2(\mathbb{R}^N)} < \varrho(m), \\
                                                               \|\nabla \gamma(1)\|^2_{L^2(\mathbb{R}^N)} >  4 \varrho(m), \\
                                                               \max \big\{ \widetilde{I}(\gamma(0)), \widetilde{I}(\gamma(1)) \big\} < \frac{1}{2} \varrho(m)
                                                            \end{array}
          \right.
        \right\}.
    \end{equation*}
  For given $m > m^{**}$, the set $\widetilde{\Gamma}_m$ is nonempty since $|\gamma| \in \widetilde{\Gamma}_m$ for any $\gamma \in \Gamma_m$. Noting that
    \begin{equation*}
      \widetilde{I}(u) = I(u^+) + \frac{1}{2} \|\nabla u^-\|^2_{L^2(\mathbb{R}^N)}
    \end{equation*}
  and repeating the proof of Lemma \ref{lemma:E_mp}, we have
    \begin{equation*}
      \widetilde{E}_{mp, m} := \inf_{\gamma \in \widetilde{\Gamma}_{m}} \sup_{t \in [0, 1]} \widetilde{I}(\gamma(t)) \geq \varrho(m) > 0.
    \end{equation*}
  Since $\widetilde{f}$ satisfies the conditions $(\hyperlink{f1}{f1})-(\hyperlink{f5}{f5})$, by an adaptation of the argument of Lemma \ref{lemma:solution}, we obtain a radial solution $w \in S_m$ of $- \Delta u + \mu u = \widetilde{f}(u)$ for some $\mu >0$. Clearly, $w^- = 0$ and thus $w \in S_m$ is a nonnegative radial solution of \eqref{problem} with $I(w) = \widetilde{I}(w) = \widetilde{E}_{mp, m} \geq \varrho(m)$.
\end{remark}

\noindent
{\bf Proof of Theorem \ref{theorem:MP}.} It follows directly from Lemmas \ref{lemma:E_mp}, \ref{lemma:solution} and Remark \ref{rmk:positive}. \hfill $\square$

%%%%%%%%%%%%%%%%%%%%%%%%%%%%%%%%%%%%%%%%%%%%%%%%%%%%%%%%%%%%%%%%%%%%%%%%%%%%%%%%%%%%%%%%%%%%%%%%%%%%%%%%%%%%%%%%%%%%%%%%%%%%%%%%%%%%%%
%%%%%%%%%%%%%%%%%%%%%%%%%%%%%%%%%%%%%%%%%%%%%%%%%%%%%%%%%%%%%%%%%%%%%%%%%%%%%%%%%%%%%%%%%%%%%%%%%%%%%%%%%%%%%%%%%%%%%%%%%%%%%%%%%%%%%%

\section{Symmetric mountain pass solutions}\label{sect:SMP}

In this section we exploit a new version of the Symmetric Mountain Pass Theorem, to prove the multiplicity result Theorem \ref{theorem:multiplicity}, which gives the existence of more and more normalized radial solutions with positive energies when the mass is getting larger. To this end, for any $m > 0$ and $k \in \mathbb{N}^+$, we define the family of odd continuous maps
  \begin{equation*}
    \Gamma_{m, k} :=
      \left\{
        \gamma \in C(\mathbb{A}_k, S_m\cap H^1_r(\mathbb{R}^N)) \left| \begin{array}{c}
                                                               \gamma[-\sigma] = - \gamma[\sigma]~\text{for any}~\sigma \in \mathbb{A}_k,\\
                                                               I(\gamma[\sigma]) < 0~\text{when}~\sigma \in \mathbb{S}^{k-1},\\
                                                               \max \big\{ \|\nabla \gamma[\sigma]\|^2_{L^2}, I(\gamma[\sigma]) \big\} < \frac{1}{2} \rho(m)~\text{if}~\sigma \in \frac{1}{2}\mathbb{S}^{k-1}
                                                            \end{array}
        \right.
      \right\},
  \end{equation*}
where $\mathbb{S}^{k-1}$ is the unit sphere in $\mathbb{R}^k$,
  \begin{equation*}
    \mathbb{A}_k := \left\{ \sigma \in \mathbb{R}^k ~|~ \frac{1}{2} \leq |\sigma| \leq 1 \right\},
  \end{equation*}
and $\rho(m) > 0$ is  given in  Lemma \ref{lemma:I} $(iii)$. Before going further, it is necessary to show the nonemptiness of $\Gamma_{m, k}$ and this seems to be a delicate issue under our mild conditions.
\begin{lemma}\label{lemma:nonemptiness}
  Assume that $N \geq 2$ and $f \in C(\mathbb{R}, \mathbb{R})$ is odd satisfying $(\hyperlink{f1}{f1})-(\hyperlink{f4}{f4})$. Then for any $k \in \mathbb{N}^+$ there exists $m_k > 0$ such that $\Gamma_{m,k} \neq \emptyset$ when $m > m_k$.
\end{lemma}
\proof Fix $k \in \mathbb{N}^+$. By  \cite[Theorem 10]{Be83-2}, there exists an odd continuous map $\pi_k : \mathbb{S}^{k-1} \to H^1_r(\mathbb{R}^N) \setminus \{0\}$ such that
  \begin{equation*}
    \inf_{\sigma \in \mathbb{S}^{k-1}} \int_{\mathbb{R}^N} F(\pi_k[\sigma])dx \geq 1.
  \end{equation*}
Then for any $m > 0$ we can define an odd continuous map $\pi_{m,k} : \mathbb{S}^{k-1} \to S_m \cap H^1_r(\mathbb{R}^N)$ as follows:
  \begin{equation*}
    \pi_{m,k}[\sigma](x) := \pi_k[\sigma]\left( m^{-1/N} \cdot \|\pi_k[\sigma]\|^{2/N}_{L^2(\mathbb{R}^N)}\cdot x \right), \qquad \sigma \in \mathbb{S}^{k-1}~\text{and}~x \in \mathbb{R}^N.
  \end{equation*}
Since $\mathbb{S}^{k-1}$ is compact and $0 \notin \pi_k[\mathbb{S}^{k-1}]$, one may find $\alpha_k, \beta_k, \beta'_k > 0$ independent of $\sigma \in \mathbb{S}^{k-1}$ such that $\|\nabla \pi_k[\sigma]\|^2_{L^2(\mathbb{R}^N)} \leq \alpha_k$ and $\beta_k \leq \|\pi_k[\sigma]\|^2_{L^2(\mathbb{R}^N)} \leq \beta'_k$. Thus
  \begin{equation*}
    \begin{split}
      I(\pi_{m,k}[\sigma])
      & = \frac{1}{2} \int_{\mathbb{R}^N} |\nabla \pi_{m, k}[\sigma]|^2 dx-\int_{\mathbb{R}^N} F(\pi_{m, k}[\sigma]) dx\\
      & = \frac{m^{\frac{N-2}{N}}}{2\|\pi_k[\sigma]\|^{2(N-2)/N}_{L^2(\mathbb{R}^N)}}\int_{\mathbb{R}^N}|\nabla \pi_k[\sigma]|^2dx - \frac{m}{\|\pi_k[\sigma]\|^2_{L^2(\mathbb{R}^N)}} \int_{\mathbb{R}^N} F(\pi_k[\sigma])dx\\
      & \leq \frac{1}{2} \alpha_k \beta^{(2-N)/N}_k \cdot m^{\frac{N-2}{N}} - (\beta'_k)^{-1} \cdot m =: g_k(m).
    \end{split}
  \end{equation*}
As a consequence, there exists $m_k > 0$ large enough such that
  \begin{equation}\label{eq:negative}
    \sup_{\sigma \in \mathbb{S}^{k-1}} I(\pi_{m,k} [\sigma]) \leq g_k(m) < 0 \qquad \text{for any}~m > m_k.
  \end{equation}
By the fact that $\|\nabla (s \diamond \pi_{m,k}[\sigma] )\|_{L^2(\mathbb{R}^N)}$ converges uniformly in $\sigma \in \mathbb{S}^{k-1}$ to zero as $s \to - \infty$, we can also choose $C_{m,k} > 0$ large enough such that
  \begin{equation}\label{eq:bound_nabla}
    \sup_{\sigma \in \mathbb{S}^{k-1}} \Bigl\|\nabla \Big( (-C_{m,k}) \diamond \pi_{m,k}[\sigma] \Big)\Bigr\|^2_{L^2(\mathbb{R}^N)} < \frac{1}{4} \rho(m)
  \end{equation}
and, using Lemma \ref{lemma:I} $(i)$,
  \begin{equation*}
    \sup_{\sigma \in \mathbb{S}^{k-1}} \left| \int_{\mathbb{R}^N} F \big( (-C_{m,k}) \diamond \pi_{m,k}[\sigma] \big) dx \right| < \frac{1}{4} \rho(m).
  \end{equation*}
In particular,
  \begin{equation}\label{eq:bound_I}
    \sup_{\sigma \in \mathbb{S}^{k-1}} I\big( (-C_{m,k}) \diamond \pi_{m,k}[\sigma] \big) < \frac{1}{2} \rho(m).
  \end{equation}
We now introduce an odd continuous map $\gamma_{m,k}: \mathbb{A}_k \to S_m \cap H^1_r(\mathbb{R}^N)$ as follows:
  \begin{equation*}
    \gamma_{m,k}[\sigma]:=\big(2C_{m,k}(|\sigma|-1)\big)\diamond \pi_{m,k}\left[\frac{\sigma}{|\sigma|}\right], \qquad \sigma \in \mathbb{A}_k.
  \end{equation*}
In view of \eqref{eq:negative}, \eqref{eq:bound_nabla} and \eqref{eq:bound_I}, it is clear that $\gamma_{m,k} \in \Gamma_{m,k}$ for any $m > m_k$.  \hfill $\square$

From now on, we fix the integer $k \in \mathbb{N}^+$ and let $m_k > 0$ be the number  given in  Lemma \ref{lemma:nonemptiness}. In order to obtain a sequence of convenient minimax values of $I_{|S_m \cap H^1_r(\mathbb{R}^N)}$, for each $j \in \{1, 2, \cdots, k\}$, we introduce the family of compact symmetric sets
  \begin{equation*}
    \Lambda_{m,j}:=\left\{\gamma[\overline{\mathbb{A}_{j+l}\setminus Y}]~\left|~\begin{array}{c}
                                                                         0\leq l\leq k-j,~\gamma\in \Gamma_{m,j+l},\\
                                                                         Y\subset \mathbb{A}_{j+l}~\text{is closed},~-Y=Y~\text{and}~\text{genus}(Y)\leq l
                                                                        \end{array}
                                                                 \right.
    \right\}.
  \end{equation*}
Here, for a given nonempty closed symmetric subset $Y$ of a Banach space, the notation $\text{genus}(Y)$ denotes the genus of $Y$ and is defined by
  \begin{equation*}
    \text{genus}(Y) := \min \big\{ n \in \mathbb{N}^+ ~|~ \exists~\varphi : Y \to \mathbb{R}^n \setminus \{0\}, ~\varphi~\text{is odd and continuous} \big\}.
  \end{equation*}
We set $\text{genus}(Y) = \infty$ if such a map $\varphi$ does not exist and set $\text{genus}(Y) = 0$ if $Y = \emptyset$.
One may refer to Sec. 7 in \cite{Ra86} for basic properties of the genus. The following properties about $\Lambda_{m,j}$ will be useful in our later arguments.
\begin{lemma}\label{lemma:Lambda_mj}
  Assume that $N \geq 2$ and $f \in C(\mathbb{R}, \mathbb{R})$ is odd satisfying $(\hyperlink{f1}{f1})-(\hyperlink{f4}{f4})$. Then for any $m > m_k$ the following statements hold.
    \begin{itemize}
      \item[$(i)$] $\Lambda_{m,1} \supset \Lambda_{m,2} \supset \cdots \supset \Lambda_{m,k} \neq \emptyset$.
      \item[$(ii)$] Let $\psi: S_m\cap H^1_r(\mathbb{R}^N) \to S_m\cap H^1_r(\mathbb{R}^N)$ be an odd continuous map such that
                      \begin{equation*}
                        \psi(u)= u \qquad \text{if}~I(u) < \frac{1}{2} \rho(m).
                      \end{equation*}
          Then $\psi(\Lambda_{m,j}) \subset \Lambda_{m,j}$ for each $j\in\{1, 2, \cdots, k\}$.
      \item[$(iii)$] If $A \in \Lambda_{m,j}$ for some $j\in\{1, 2, \cdots, k\}$ and $B \subset S_m \cap H^1_r(\mathbb{R}^N)$ is a closed symmetric set with $\text{genus}(B) \leq i <j$, then
                       \begin{equation*}
                         \overline{A \setminus B} \in \Lambda_{m,j-i}.
                       \end{equation*}
      \item[$(iv)$] For any $j\in\{1, 2, \cdots, k\}$ and $A \in \Lambda_{m,j}$ there exists $w \in A$ such that
                      \begin{equation*}
                        \|\nabla w\|^2_{L^2(\mathbb{R}^N)} = 4\rho(m).
                      \end{equation*}
    \end{itemize}
\end{lemma}
\proof $(i)$ It follows from the definition of $\Lambda_{m,j}$ and Lemma \ref{lemma:nonemptiness}.

$(ii)$  We only need to prove that $\psi(A) \in \Lambda_{m,j}$ for any $j\in\{1, 2, \cdots, k\}$ and $A \in \Lambda_{m,j}$. Denote $A = \gamma[\overline{\mathbb{A}_{j+l}\setminus Y}] \in \Lambda_{m,j}$ for some $\gamma\in \Gamma_{m,j+l}$ and $Y \subset \mathbb{A}_{j+l}$. It is clear that
  \begin{equation*}
    \psi \circ \gamma \in \Gamma_{m,j+l}
  \end{equation*}
and thus $\psi(A) =  \psi \circ \gamma [\overline{\mathbb{A}_{j+l}\setminus Y}] \in \Lambda_{m,j}$.

$(iii)$ We adapt the argument in \cite[Proposition 9.18]{Ra86}. Suppose $A = \gamma[\overline{\mathbb{A}_{j+l}\setminus Y}] \in \Lambda_{m,j}$,  where $\gamma\in \Gamma_{m,j+l}$ and $Y \subset \mathbb{A}_{j+l}$ with $\text{genus}(Y)\leq l \leq k-j$. Let $B$ be as above with $\text{genus}(B) \leq i <j$. We claim that
  \begin{equation}\label{eq:claim}
    \overline{A \setminus B} = \gamma \big[\overline{\mathbb{A}_{j+l} \setminus (Y \cup \gamma^{-1}(B) )}\big].
  \end{equation}
Indeed, if $u \in \gamma \big[\mathbb{A}_{j+l} \setminus (Y \cup \gamma^{-1}(B)) \big]$, then
  \begin{equation*}
    u \in \gamma[\mathbb{A}_{j+l} \setminus Y] \setminus B \subset A\setminus B \subset \overline{A \setminus B}.
  \end{equation*}
Therefore
  \begin{equation}\label{eq:claim_1}
    \gamma \big[\mathbb{A}_{j+l} \setminus (Y \cup \gamma^{-1}(B)) \big] \subset \overline{A \setminus B}.
  \end{equation}
On the other hand, for any $u \in A \setminus B$ there exists
  \begin{equation*}
    \sigma = \sigma_u \in \overline{ \mathbb{A}_{j+l} \setminus Y} \setminus \gamma^{-1}(B) \subset \overline{\mathbb{A}_{j+l} \setminus (Y \cup \gamma^{-1}(B) )}
  \end{equation*}
such that $u = \gamma[\sigma]$. Thus
  \begin{equation}\label{eq:claim_2}
    A \setminus B \subset \gamma \big[\overline{\mathbb{A}_{j+l} \setminus (Y \cup \gamma^{-1}(B) )}\big].
  \end{equation}
Noting that $\gamma$ is continuous, the claim \eqref{eq:claim} follows from \eqref{eq:claim_1} and \eqref{eq:claim_2}. Since $\gamma$ is an odd continuous map, $\gamma^{-1}(B)$ and then $Y \cup \gamma^{-1}(B)$ are both closed and symmetric. By some basic properties of the genus,
  \begin{equation*}
    \begin{split}
      \text{genus}(Y \cup \gamma^{-1}(B))
        & \leq \text{genus}(Y) + \text{genus}(\gamma^{-1}(B)) \\
        & \leq \text{genus}(Y) + \text{genus}(B) \\
        & \leq  l + i.
    \end{split}
  \end{equation*}
In view of \eqref{eq:claim}, we finally conclude that $\overline{A \setminus B} \in \Lambda_{m,j-i}$.

$(iv)$ The proof is motivated by that of \cite[Proposition 9.23]{Ra86}. Set $A = \gamma[\overline{\mathbb{A}_{j+l}\setminus Y}] \in \Lambda_{m,j}$,  where $\gamma\in \Gamma_{m,j+l}$ and $Y \subset \mathbb{A}_{j+l}$ with $\text{genus}(Y)\leq l$. For any $\sigma \in \mathbb{S}^{j + l -1}$, by the fact that $I(\gamma[\sigma]) < 0$ and Lemma \ref{lemma:I} $(iii)$, we have $\|\nabla \gamma[\sigma]\|^2_{L^2(\mathbb{R}^N)} > 4\rho(m)$. Denote
  \begin{equation*}
    \Omega_m:= \big \{u \in S_m \cap H^1_r(\mathbb{R}^N)~|~\|\nabla u\|^2_{L^2(\mathbb{R}^N)} < 4\rho(m) \big \}
  \end{equation*}
and let $U \subset \mathbb{A}_{j + l}$ be the connected component of $\gamma^{-1}(\Omega_m)$ containing $\frac{1}{2} \mathbb{S}^{j+l-1}$. Noting that $\mathcal{O} := U \cup \{ \sigma \in \mathbb{R}^{j+l}~|~|\sigma| < 1/2 \}$ is a bounded symmetric neighborhood of $0 \in \mathbb{R}^{j+l}$, we have $\text{genus}(\partial \mathcal{O}) = j + l$ and thus
  \begin{equation*}
    \text{genus}(\overline{\partial \mathcal{O} \setminus Y}) \geq \text{genus}(\partial \mathcal{O}) - \text{genus}(Y) \geq j \geq 1.
  \end{equation*}
In particular, $\overline{\partial \mathcal{O} \setminus Y} \neq \emptyset$. Since
  \begin{equation*}
    \gamma[\overline{\partial \mathcal{O} \setminus Y}] \subset A \cap \partial \Omega_m,
  \end{equation*}
there exists  $w \in A$ such that $\|\nabla w\|^2_{L^2(\mathbb{R}^N)} = 4\rho(m)$. \hfill $\square$

\begin{lemma}\label{lemma:E_mj}
  Assume that $N \geq 2$ and $f \in C(\mathbb{R}, \mathbb{R})$ is odd satisfying $(\hyperlink{f1}{f1})-(\hyperlink{f4}{f4})$. For any $m > m_k$ we define minimax values
    \begin{equation*}
      E_{m, j} := \inf_{A \in \Lambda_{m,j}}\max_{u \in A} I(u), \qquad  j = 1, 2, \cdots, k.
    \end{equation*}
  Then $0 < \rho(m) \leq E_{m,1} \leq E_{m,2}\leq \cdots \leq E_{m,k} < + \infty$.
\end{lemma}
\proof By Lemma \ref{lemma:Lambda_mj} $(i)$, we have $E_{m,1} \leq E_{m,2}\leq \cdots \leq E_{m,k} < +\infty$. From Lemma \ref{lemma:Lambda_mj} $(iv)$ and  Lemma \ref{lemma:I} $(iii)$, it follows that
  \begin{equation*}
    \max_{u \in A} I(u) \geq \rho(m) \qquad \text{for any}~A \in \Lambda_{m,j}
  \end{equation*}
and thus $E_{m,j} \geq \rho(m)$. \hfill $\square$

As a final preparation for the proof of Theorem \ref{theorem:multiplicity}, using essentially the deformation result Lemma \ref{lemma:deformation} and the compactness result Lemma \ref{lemma:PSP_c}, we manage to prove  a key lemma which not only shows that $E_{m, j}$ is indeed a critical value for any $j=1,2, \cdots, k$ but also makes an appropriate multiplicity statements about degenerate critical values.
\begin{lemma}\label{lemma:solutions}
  Assume that $N \geq 2$ and $f \in C(\mathbb{R}, \mathbb{R})$ is odd satisfying $(\hyperlink{f1}{f1})-(\hyperlink{f5}{f5})$. Let $m > m_k$ and recall that $K^c$ stands for the set of  critical points of $I_{|S_m\cap H^1_r(\mathbb{R}^N)}$  at a level $c \in \mathbb{R}$. If
    \begin{equation*}
      E_{m,j} = E_{m,j+1} = \cdots = E_{m,j+i} =: c \geq \rho(m)
    \end{equation*}
  for some $j \geq 1$ and $i \geq 0$ satisfying $i + j \leq k$, then
    \begin{equation*}
      \text{genus}(K^c) \geq i + 1.
    \end{equation*}
  In particular, $I_{|S_m \cap H^1_r(\mathbb{R}^N)}$ has infinitely many distinct critical points at the level $c$ if $i \geq 1$.
\end{lemma}
\proof Assume by contradiction that $\text{genus}(K^c) \leq i$. Since $K^c$ is compact by Lemma \ref{lemma:PSP_c}, we have a neighborhood $\mathcal{N} \subset S_m \cap H^1_r(\mathbb{R}^N)$ of $K^c$ such that $\text{genus}(\mathcal{N}) = \text{genus}(K^c) \leq i$. Applying the deformation result Lemma \ref{lemma:deformation} with $\mathcal{O} = \mathcal{N}$ and $\overline{\varepsilon} = \rho(m)/2 > 0$, there exists $\varepsilon \in (0, \overline{\varepsilon})$ and $\eta \in C([0,1] \times (S_m \cap H^1_r(\mathbb{R}^N)), S_m \cap H^1_r(\mathbb{R}^N))$ such that $\eta(1, \cdot)$ is odd and
  \begin{equation}\label{eq:deformation_SMP}
    \eta(1, I^{c + \varepsilon} \setminus \mathcal{O}) \subset I^{c - \varepsilon}.
  \end{equation}
By the definition of $E_{m, j + i}$, we choose $A \in \Lambda_{m, j + i}$ such that
  \begin{equation}\label{eq:A}
    \max_{u \in A} I(u) \leq c + \varepsilon.
  \end{equation}
From Lemma \ref{lemma:Lambda_mj} $(iii)$, $\overline{A \setminus \mathcal{O}} \in \Lambda_{m, j}$. For any $u \in S_m \cap H^1_r(\mathbb{R}^N)$ satisfying $I(u) < \rho(m)/2$, we have
  \begin{equation*}
    I(u) < \frac{1}{2}\rho(m) \leq c - \overline{\varepsilon},
  \end{equation*}
and thus $\eta(1, u) = u$ by Item $(ii)$ of Lemma \ref{lemma:deformation}. In view of Lemma \ref{lemma:Lambda_mj} $(ii)$, it follows that $\eta(1, \overline{A \setminus \mathcal{O}}) \in \Lambda_{m, j}$. Now, by the definition of $E_{m, j}$, \eqref{eq:deformation_SMP} and \eqref{eq:A}, we obtain
  \begin{equation*}
    c = E_{m, j} \leq \max_{u \in \eta(1, \overline{A \setminus \mathcal{O}})} I(u) \leq c - \varepsilon
  \end{equation*}
which is a contradiction. \hfill $\square$

With all the technical results in place, we can now prove Theorem \ref{theorem:multiplicity}.

\medskip
\noindent
{\bf Proof of Theorem \ref{theorem:multiplicity}.} For any $k \in \mathbb{N}^+$, let $m_k > 0$ be the number given in  Lemma \ref{lemma:nonemptiness}. Then the conclusion follows from Lemmas \ref{lemma:E_mj} and \ref{lemma:solutions}. \hfill $\square$

\medskip
We end this section by observing that by an appropriate adaptation of the arguments for Theorem \ref{theorem:multiplicity}, a nonradial variant of that multiplicity result can be established when $N = 4$ or $N \geq 6$.
\begin{theorem}[Nonradial sign-changing solutions at positive levels]\label{theorem:nonradial}
  Assume $N =4$ or $N \geq 6$, and $f \in C(\mathbb{R}, \mathbb{R})$ is odd satisfying $(\hyperlink{f1}{f1})-(\hyperlink{f5}{f5})$. Then for each $k \in \mathbb{N}^+$ there exists $m_k>0$ such that when $m > m_k$ the problem \eqref{problem} has at least $k$ distinct nonradial sign-changing solutions with positive energies.
\end{theorem}
One should note here that the working space for proving Theorem \ref{theorem:nonradial} is the same as the one used in \cite[Theorem 1.2]{JL19} and that the special map in \cite[Lemma 3.4]{JL19} can be used to show the nonemptiness of a certain substitute to the family of odd continuous maps $\Gamma_{m, k}$. The details of the proof are left to the interested reader.

%%%%%%%%%%%%%%%%%%%%%%%%%%%%%%%%%%%%%%%%%%%%%%%%%%%%%%%%%%%%%%%%%%%%%%%%%%%%%%%%%%%%%%%%%%%%%%%%%%%%%%%%%%%%%%%%%%%%%%%%%%%%%%%%%%%%%%
%%%%%%%%%%%%%%%%%%%%%%%%%%%%%%%%%%%%%%%%%%%%%%%%%%%%%%%%%%%%%%%%%%%%%%%%%%%%%%%%%%%%%%%%%%%%%%%%%%%%%%%%%%%%%%%%%%%%%%%%%%%%%%%%%%%%%%

\section{Application to the cubic-quintic nonlinear Schr\"{o}dinger equation}\label{sect:application}

In this last section we present some implications of the results of the preceding sections to the three-dimensional cubic-quintic nonlinear Schr\"{o}dinger equation
  \begin{equation}\label{eq:cubic-quintic-evolution}
    i \partial_t \psi = - \Delta \psi - |\psi|^2 \psi + |\psi|^4 \psi, \qquad (t,x) \in \mathbb{R} \times \mathbb{R}^3,
  \end{equation}
subject to the initial data $\psi_{| t=0} = \psi_0 \in H^1(\mathbb{R}^3)$. The cubic-quintic nonlinearity was introduced in \cite{PPT79} and is now used in several physical models. We refer to the review \cite{M19} for precise references and also to the papers \cite{CKS21,CS21,KMV21,KOPV17,LRN20, TVZ07} for recent results.  Note that, although  \eqref{eq:cubic-quintic-evolution} is energy critical, the global existence of solutions holds for any initial data,  see \cite[Theorem 1.1]{KOPV17} or \cite{Z06}.

When looking to standing waves to \eqref{eq:cubic-quintic-evolution}, namely solutions of the form $\psi(t, x) = e^{i \omega t}u(x)$ with $\omega \in \mathbb{R}$ and $u \in H^1(\mathbb{R}^3)$, we are led  to the equation
  \begin{equation}\label{eq:cubic-quintic-static}
    - \Delta u + \omega u =  |u|^2 u - |u|^4 u \quad \mbox{in} ~ H^1(\mathbb{R}^3).
  \end{equation}
The equation \eqref{eq:cubic-quintic-static} is of the form of \eqref{eq:equation-libre} with the choice $N = 3$ and $f(t) = |t|^{2}t - |t|^{4}t$. It is readily seen that the conditions $(\hyperlink{f1}{f1}) - (\hyperlink{f5}{f5})$ and \eqref{eq:key_3} are satisfied, and thus all the results established in the previous sections apply.

In order to discuss more details, for fixed $\omega > 0$  we denote by $J_\omega$ the $C^1$ action functional
  \begin{equation*}
    J_\omega(u)
      := I(u) + \frac{\omega}{2} \|u\|^2_{L^2(\mathbb{R}^3)}
      = \frac{1}{2}\| \nabla u\|^2_{L^2(\mathbb{R}^3)} + \frac{\omega}{2} \|u\|^2_{L^2(\mathbb{R}^3)} - \frac{1}{4} \|u\|^4_{L^4(\mathbb{R}^3)} + \frac{1}{6} \|u\|^6_{L^6(\mathbb{R}^3)}.
  \end{equation*}
A nontrivial solution $v \in H^1(\mathbb{R}^3) \setminus \{0\}$ to \eqref{eq:cubic-quintic-static} is called an action ground state if $J_\omega(v) \leq J_\omega(u)$ for any nontrivial solution $u$ to \eqref{eq:cubic-quintic-static}. In \cite[Theorem 2.2]{KOPV17} the following basic properties were established, see also \cite{LRN20}.

\begin{proposition}\label{prop: results KOP}
For any given frequency $\omega \in (0, 3/16)$, there exists a unique positive radial solution $U_\omega \in H^1(\mathbb{R}^3)$ to \eqref{eq:cubic-quintic-static}. Moreover,
  \begin{itemize}
    \item[$(i)$] the map $\omega \mapsto U_\omega$ is real analytic;
    \item[$(ii)$] $U_\omega$ is an action ground state to \eqref{eq:cubic-quintic-static}.
  \end{itemize}
\end{proposition}
Many complementary results have also been obtained.  In particular, the precise asymptotic behavior of $U_\omega$ in the limits $\omega \to 0$ and $\omega \to 3/16$, see \cite[Theorem 2.2]{KOPV17} and  \cite[Theorems 3 and 4]{LRN20}.  From these asymptotics, it is proved that $U_\omega$ is  orbitally unstable and orbitally stable when $\omega$ is sufficiently close to $0$ and $3/16$ respectively, see \cite[Theorems 3 and 4]{LRN20}.

\begin{figure}[htbp]
  \centering
  \includegraphics[width=0.75\textwidth]{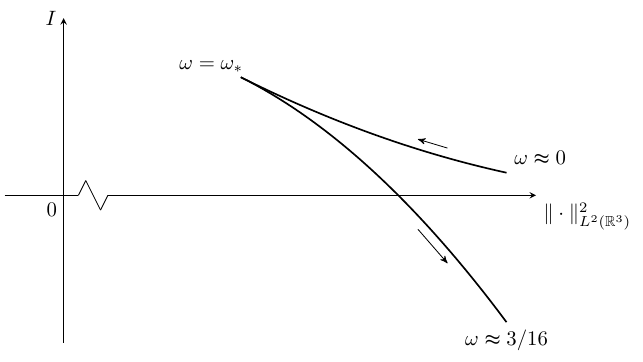}
  \caption{Schematic depiction of the mass-energy curve $(\|U_\omega\|^2_{L^2(\mathbb{R}^3)}, I(U_\omega))$ for the action ground state $U_\omega$, based on numerics.}
  \label{fig:E-M}
\end{figure}

However, open questions remain.  Among other things, in view of several numerical simulations in \cite{CKS21, DMM00, KOPV17,LRN20, MMCML00}, see in particular \cite[Fig. 2]{KOPV17} and \cite[Fig. 3]{LRN20}, two main conjectures have been formulated.
  \begin{itemize}
    \item[(C.1)] \hypertarget{C.1}{} There is an $\omega_* \in (0, 3/16)$ such that the mass function $\omega \mapsto \|U_\omega\|^2_{L^2(\mathbb{R}^3)}$ is strictly decreasing for $\omega \in (0, \omega_*)$ and strictly increasing for $\omega \in (\omega_*, 3/16)$.
    \item[(C.2)] \hypertarget{C.2}{} The behavior of the mass-energy curve $(\|U_\omega\|^2_{L^2(\mathbb{R}^3)}, I(U_\omega))$ is of the form of Figure \ref{fig:E-M}.
  \end{itemize}

\begin{remark}
  If the conjecture \emph{(\hyperlink{C.1}{C.1})} proves correct then,  as a rather direct consequence of classical results of \cite{GSS87}, one can deduce that the standing waves associated to $U_\omega$ are orbitally unstable if $\omega \in (0, \omega_*)$ and orbitally stable when $\omega \in (\omega_*, 3/16)$. In that direction, note that for stability issues, even for the existence of action ground states, the equations \eqref{eq:cubic-quintic-evolution} and \eqref{eq:cubic-quintic-static} have actually to be considered in $H^1(\mathbb{R}^3, \mathbb{C})$, but  the discussion can be reduced to $H^1(\mathbb{R}^3, \mathbb{R})$ by standard arguments.
\end{remark}

Combining Proposition \ref{prop: results KOP} and the results of the previous sections we obtain Theorem \ref{theorem:cubic-quintic} below.  With respect to the preceding studies,  it provides additional information on the action ground states $U_\omega$ in ranges of the frequency $\omega$ in which $U_\omega$ is not a global minimizer of the constrained energy functional. In particular, it appears that $U_\omega$ can be a local minimizer or of saddle type. Also, new information on the sign of $I(U_\omega)$ is obtained.

\begin{theorem}\label{theorem:cubic-quintic}
  Let $N=3$, $f(t) = |t|^2t - |t|^4 t$, and $E_m \leq 0$, $m^* > 0$, $\overline{E}_m \geq 0$ and $m^{**} \in (0, m^*)$ be the numbers given in Theorems \ref{th-1} and \ref{theorem:local}. Then the following statements hold.
    \begin{itemize}
      \item[$(i)$] For any $m \in (m^{**}, m^*)$, the set of minimizers of the local infimum $\overline{E}_m > 0$ is nonempty and orbitally stable.
      \item[$(ii)$] For any $m > m^{**}$, there exist two distinct frequencies $\omega_1, \omega_2 \in (0, 3/16)$ such that the associated free problems \eqref{eq:cubic-quintic-static} have respectively the action ground states $v_m$ and $w_m$ such that $\|v_m\|^2_{L^2(\mathbb{R}^3)} = \|w_m\|^2_{L^2(\mathbb{R}^3)} = m$ and
            \begin{equation*}
              \left\{
                     \begin{aligned}
                       & I(w_m) > 0 > E_m =I(v_m), & \qquad  & \text{for} ~ m > m^*, \\
                       & I(w_m) > 0 = \overline{E}_m = I(v_m), & \qquad & \text{for} ~ m = m^*,\\
                       & I(w_m) > I(v_m) = \overline{E}_m > 0, & \qquad & \text{for} ~ m \in (m^{**}, m^*).
                     \end{aligned}
              \right.
            \end{equation*}
          In particular, $v_m$ and $w_m$ correspond to $U_{\omega_1}$ and $U_{\omega_2}$ respectively.
      \item[$(iii)$] For any $m \in (m^{**}, m^*)$, the action ground state $v_m = U_{\omega_1}$ is a local rather than a global minimizer of the constrained energy function $I_{|S_m}$.
      \item[$(iv)$] For any $m > m^{**}$, the action ground state $w_m = U_{\omega_2} \in S_m \cap H^1_r(\mathbb{R}^3)$ corresponds to a mountain pass level of the constrained energy function $I_{|S_m \cap H^1_r(\mathbb{R}^3)}$.
    \end{itemize}
\end{theorem}
\proof For any $m \in (m^{**}, m^*)$, by Theorem \ref{theorem:local} and Remark \ref{rmk:comments} $(i)$, the local infimum $\overline{E}_m$ is achieved and any associated minimizing sequence is, up to a subsequence and up to translations in $\mathbb{R}^3$, strongly convergent. Following the strategy laid down in \cite{CL82} and recalling the fact that the global existence of solutions to \eqref{eq:cubic-quintic-evolution} holds for any initial data, it is not difficult to see that the set of local minimizers is orbitally stable. This proves Item $(i)$. Since $f(t) = |t|^{2}t - |t|^{4}t$ is odd,  the normalized solutions $v \in S_m$ and $w \in S_m \cap H^1_r(\mathbb{R}^3)$ obtained in Theorems \ref{th-1}, \ref{theorem:local} and \ref{theorem:MP} can be assumed to be nonnegative. Noting that the associated Lagrange multipliers $\mu(m,v)$ and $\mu(m,w)$ are positive, by regularity and the strong maximum principle, we have that they are strictly positive on $\mathbb{R}^3$. In view of the well known symmetry result in \cite{GNN81}, the positive solution $v$ can be assumed further to be radially symmetric. Now, with the choice $(\omega_1, v_m) := (\mu(m, v), v)$ and $(\omega_2, w_m) := (\mu(m, w), w)$, it is clear that Proposition \ref{prop: results KOP} and Theorem \ref{theorem:MP} imply Item $(ii)$. Finally, Items $(iii)$ and $(iv)$ follows directly from the proof of Item $(ii)$. \hfill $\square$

\begin{remark}
  \begin{itemize}
    \item[$(i)$] Concerning Item $(i)$ in Theorem \ref{theorem:cubic-quintic},  note that in \cite[Theorem 5.6]{KOPV17} the existence of such a piece of branch was already obtained but the issue of stability was not considered. We also mention that in \cite[Theorem 5.2]{KOPV17}, some properties of the map $m \mapsto \overline{E}_m$ presented in our Theorem \ref{theorem:local} were derived for this specific nonlinearity.
    \item[$(ii)$] Theorem \ref{theorem:cubic-quintic} $(ii)$ can be considered as a somewhat new evidence for the conjectures \emph{(\hyperlink{C.1}{C.1})} and \emph{(\hyperlink{C.2}{C.2})} to be true, even though little is known about $\omega_1, \omega_2 \in (0, 3/16)$.
    \item[$(iii)$] Items $(iii)$ and $(iv)$ of Theorem \ref{theorem:cubic-quintic} demonstrate that for some frequency $\omega \in (0, 3/16)$ the action ground state $U_\omega$ is not a minimizer of the global infimum energy $E_m$ with $m = \|U_\omega\|^2_{L^2(\mathbb{R}^3)}$. This is consistent with \cite[Theorem 2.5]{CKS21} and in particular we present new counterexamples that disprove the converse statement of \cite[Theorem 1.6 $(i)$]{JL21}.
    \item[$(iv)$] In view of Theorem \ref{theorem:cubic-quintic} $(i)$, it seems reasonable to expect that when $m \in (m^{**}, m^*)$ the standing wave $e^{i \omega_1 t}U_{\omega_1}$ is orbitally stable. On the other hand, for any $m > m^*$, the variational characterization of $w_m = U_{\omega_2}$ inclines to believe that $e^{i \omega_2 t}U_{\omega_2}$ is orbitally unstable. If this is indeed the case the instability will not occur by finite time blow-up as a consequence of the global existence of solutions, for any initial data, of the evolution problem \eqref{eq:cubic-quintic-evolution}. We refer to  \cite{CKS21} for a numerical exploration of this instability.
  \end{itemize}
\end{remark}

\begin{remark}\label{rmk:no small solution}
In \cite[Theorem 5.2]{KOPV17} it is proved that \eqref{eq:cubic-quintic-static} has no solution with small $L^2$ norm, showing that, in general,  a normalized solution  cannot be expected for arbitrarily small value of $m > 0$ in our Theorem \ref{theorem:local}.  Also, in \cite[Theorem 5.6]{KOPV17}, it was observed that the existence of an element $u \in S_m$ satisfying $P(u) = 0$ is not sufficient to guarantee the existence of an energy ground state solution to \eqref{problem}.  This is another indication that deriving a sharp estimate on the value $m^{**} >0$ which appears in Theorem \ref{theorem:local} may prove challenging.
\end{remark}

\begin{remark}\label{rmk:2D case}
  The two dimensional counterpart of \eqref{eq:cubic-quintic-evolution} also has  a physical interest and it was considered in the works \cite{CKS21,CS21}. Clearly, when $N=2$, the nonlinearity $f(t) = |t|^{2}t - |t|^{4}t$  satisfies $(\hyperlink{f1}{f1}) - (\hyperlink{f3}{f3})$  and \eqref{eq:key_2}, and thus the conclusions of \cite[Theorems 1.4 and 1.6]{JL21} as well as of Theorem \ref{th-1} apply. However, it does not satisfy $(\hyperlink{f4}{f4})$.  One may observe  that if $u \in S_m$ is a constrained critical point then
    \begin{equation*}
      I(u) = I(u) - \frac{1}{2}P(u) =  - \frac{1}{6} \int_{\mathbb{R}^2}|u|^6 dx  <0.
    \end{equation*}
  Thus there is no critical point at a positive energy level. Also $E_{m^*} = 0$ is not achieved.
\end{remark}

%%%%%%%%%%%%%%%%%%%%%%%%%%%%%%%%%%%%%%%%%%%%%%%%%%%%%%%%%%%%%%%%%%%%%%%%%%%%%%%%%%%%%%%%%%%%%%%%%%%%%%%%%%%%%%%%%%%%%%%%%%%%%%%%%%%%%%
%%%%%%%%%%%%%%%%%%%%%%%%%%%%%%%%%%%%%%%%%%%%%%%%%%%%%%%%%%%%%%%%%%%%%%%%%%%%%%%%%%%%%%%%%%%%%%%%%%%%%%%%%%%%%%%%%%%%%%%%%%%%%%%%%%%%%%

\section*{Acknowledgements}
\addcontentsline{toc}{section}{Acknowledgements}

Sheng-Sen Lu thanks Dr. Hao Luo at Peking University for Figure \ref{fig:E-M}. He also acknowledges the support of the China Postdoctoral Science Foundation (No. 2020M680174) and the National Natural Science Foundation of China (Nos. 11771324 and 11831009).

%%%%%%%%%%%%%%%%%%%%%%%%%%%%%%%%%%%%%%%%%%%%%%%%%%%%%%%%%%%%%%%%%%%%%%%%%%%%%%%%%%%%%%%%%%%%%%%%%%%%%%%%%%%%%%%%%%%%%%%%%%%%%%%%%%%%%%
%%%%%%%%%%%%%%%%%%%%%%%%%%%%%%%%%%%%%%%%%%%%%%%%%%%%%%%%%%%%%%%%%%%%%%%%%%%%%%%%%%%%%%%%%%%%%%%%%%%%%%%%%%%%%%%%%%%%%%%%%%%%%%%%%%%%%%

\end{document}